\documentclass{amsart}

\usepackage{amsmath}
\usepackage[latin1]{inputenc}
\usepackage{amsfonts}
\usepackage{amssymb}
\usepackage{graphicx,epsfig}
\usepackage{amsthm}
\usepackage{amscd}
\usepackage[all]{xy}
\usepackage[latin1]{inputenc}

\usepackage[T1]{fontenc}

\usepackage{amsmath}

\usepackage{amsfonts}

\usepackage{amssymb}

\usepackage{amsthm}

\usepackage{graphics}

\newcommand{\mk}{\medskip}

\newcommand{\ZZ}{\mathbb{Z}}

\newcommand{\CC}{\mathbb{C}}

\newcommand{\NN}{\mathbb{N}}

\newcommand{\QQ}{\mathbb{Q}}

\newcommand{\Glie}{\mathfrak{g}}

\newcommand{\Yim}{\mathcal{Y}}

\newcommand{\Hlie}{\mathfrak{h}}

\newcommand{\demo}{\noindent {\it \small Proof:}\quad}

\renewcommand{\NN}{\ensuremath{\mathbb{N}}}

\renewcommand{\CC}{\ensuremath{\mathbb{C}}}

\renewcommand{\QQ}{\ensuremath{\mathbb{Q}}}

\newcommand{\U}{\mathcal{U}}

\newtheorem{thm}{Theorem}[section]

\newtheorem{defi}[thm]{Definition}

\newtheorem{cor}[thm]{Corollary}

\newtheorem{prop}[thm]{Proposition}

\newtheorem{lem}[thm]{Lemma}

\author{David Hernandez}
\address{David Hernandez: \'Ecole Normale Sup\'erieure - DMA, 45, Rue d'Ulm F-75230 PARIS,
  Cedex 05  FRANCE}
\email{David.Hernandez@ens.fr\\ URL: http://www.dma.ens.fr/$\sim$dhernand}

\usepackage[all]{xy}
\title{Representations of Quantum Affinizations and Fusion Product}

\begin{document}

\begin{abstract} In this paper we study general quantum affinizations $\U_q(\hat{\Glie})$ of symmetrizable quantum Kac-Moody algebras and we develop their representation theory. We prove a triangular decomposition and we give a classication of (type $1$) highest weight simple integrable representations analog to Drinfel'd-Chari-Pressley one. A generalization of the $q$-characters morphism, introduced by Frenkel-Reshetikhin for quantum affine algebras, appears to be a powerful tool for this investigation. For a large class of quantum affinizations (including quantum affine algebras and quantum toroidal algebras), the combinatorics of $q$-characters give a ring structure $*$ on the Grothendieck group $\text{Rep}(\U_q(\hat{\Glie}))$ of the integrable representations that we classified. We propose a new construction of tensor products in a larger category by using the Drinfel'd new coproduct (it can not directly be used for $\text{Rep}(\U_q(\hat{\Glie}))$ because it involves infinite sums). In particular we prove that $*$ is a fusion product (a product of representations is a representation).\end{abstract}

\maketitle

\tableofcontents

\section{Introduction} 

In this paper $q\in\CC^*$ is not a root of unity.

\noindent V.G. Drinfel'd \cite{Dri1} and M. Jimbo \cite{jim} associated, independently, to any symmetrizable Kac-Moody algebra $\Glie$ and $q\in\CC^*$ a Hopf algebra $\U_q(\Glie)$ called quantum Kac-Moody algebra. The structure of the Grothendieck ring of integrable representations is well understood : it is analogous to the classical case $q=1$.

\noindent The quantum algebras of finite type $\U_q(\Glie)$ ($\Glie$ of finite type) have been intensively studied (see for example \cite{Cha2, lu, ro} and references therein). The quantum affine algebras $\U_q(\hat{\Glie})$ ($\hat{\Glie}$ affine algebra) are also of particular interest : they have two realizations, the usual Drinfel'd-Jimbo realization and a new realization (see \cite{Dri2, bec}) as a quantum affinization of a quantum algebra of finite type $\U_q(\Glie)$. The finite dimensional representations of quantum affine algebras are the subject of intense research (see among others \cite{Aka, Cha0, Cha, Cha2, em, Fre, Fre2, Naams, Nab, var2} and references therein). In particular they were classified by Chari-Pressley \cite{Cha, Cha2}, and Frenkel-Reshetikhin \cite{Fre} introduced the $q$-characters morphism which is a powerful tool for the study of these representations (see also \cite{kn, Fre2}).

\noindent The quantum affinization process (that Drinfel'd \cite{Dri2} described for constructing the second realization of a quantum affine algebra) can be extended to all symmetrizable quantum Kac-Moody algebras $\U_q(\Glie)$ (see \cite{jin, Naams}). One obtains a new class of algebras called quantum affinizations : the quantum affinization of $\U_q(\Glie)$ is denoted by $\U_q(\hat{\Glie})$. The quantum affine algebras are the simplest examples and are very special because they are also quantum Kac-Moody algebras. When $C$ is affine, the quantum affinization $\U_q(\hat{\Glie})$ is called a quantum toroidal algebra. It is known not to be a quantum Kac-Moody algebra but it is also of particular interest (see for example \cite{gkv, mi, mi2, Naams, Nad, Sa, Sc, stu, tu, var} and references therein). This setting is summed up in this picture :

\begin{center}
\epsfig{file=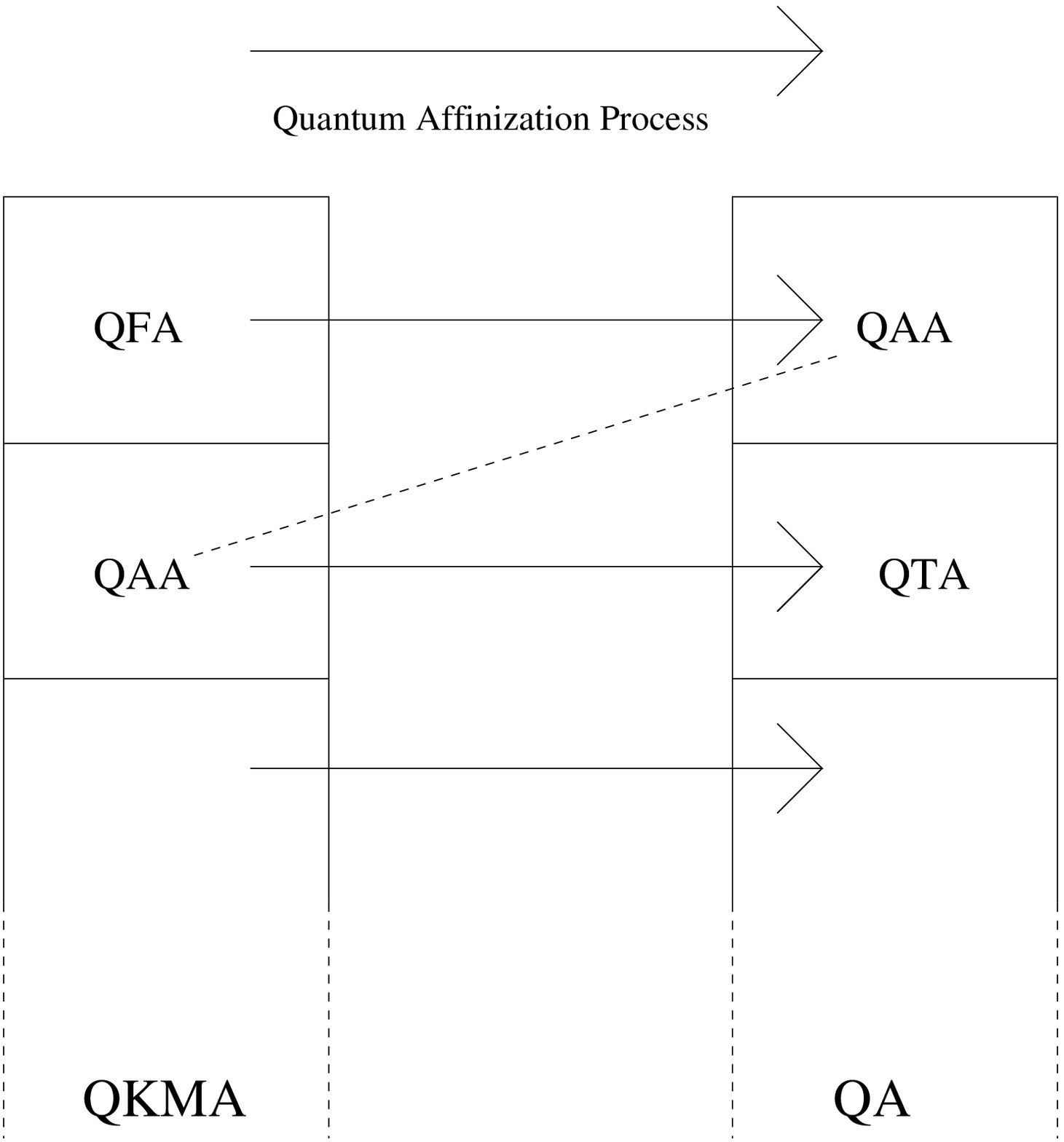,width=.6\linewidth}
\end{center}

\noindent (QKMA : Quantum Kac-Moody Algebras, QFA : Quantum Algebras of Finite type, QAA : Quantum Affine Algebras, QTA : Quantum Toroidal Algebras, QA : Quantum Affinizations; the line between the two QAA symbolizes the Drinfel'd-Beck correspondence.)

\noindent In \cite{Naams} Nakajima gave a classification of (type $1$) simple integrable highest weight modules of $\U_q(\hat{\Glie})$ when $\Glie$ is symmetric. The case $C$ of type $A_n^{(1)}$ (toroidal $\hat{sl_n}$-case) was also studied by Miki in \cite{mi} (a coproduct is also used with an approach specific to the $A_n^{(1)}$-case; but it is technically different from the general construction proposed in this paper). In \cite{her03} we proposed a combinatorial construction of q-characters (and also of their t-deformations) for generalized Cartan matrix $C$ such that $i\neq j\Rightarrow C_{i,j}C_{j,i}\leq 3$ (it includes finite and affine types except $A_1^{(1)}$, $A_2^{(2)}$); we conjectured that they were linked with a general representation theory. But in general little is known about the representation theory outside the case of quantum affine algebras.

In this paper we study general quantum affinizations and we develop their representation theory. First we prove a triangular decomposition of $\U_q(\hat{\Glie})$. We classify the (type $1$) simple highest weight integrable representations, we define and study a generalization of the morphism of $q$-characters $\chi_q$ which appears to be a natural tool for this investigation (the approach is different from \cite{her03} because $q$-characters are obtained from the representation theory and not from purely combinatorial constructions). If the quantized Cartan matrix $C(z)$ is invertible (it includes all quantum affine algebras and quantum toroidal algebras), a symmetry property of those $q$-characters with respect to the action of screening operators is proved (analog of the invariance for the action of the Weyl group in classical finite cases; the result is proved in \cite{Fre2} for quantum affine algebras); in particular those $q$-characters are the combinatorial objects considered in \cite{her03}. Moreover we get that the image of $\chi_q$ is a ring and we can define a formal ring structure on the Grothendieck group. Although quantum affine algebras are Hopf algebras, in general no coproduct has been defined for quantum affinizations (this point was also raised by Nakajima in \cite{Nad}). Drinfel'd gave formulas for a new coproduct which can be written for all quantum affinizations. They can not directly be used to define a tensor product of representations because they involve infinite sums. We propose a new construction of tensor products in a larger category with a generalization of the new Drinfel'd coproduct. We define a specialization process which allows us to interpret the ring structure that we defined on the Grothendieck group : we prove that it is a fusion product, that is to say that a product of representations is a representation (see \cite{Fu} for generalities on fusion rings and physical motivations).

In more details, this paper is organized as follows : 

\noindent in section \ref{bck} we recall backgrounds on quantum Kac-Moody algebras. In section \ref{genaff} we recall the definition of quantum affinizations and we prove a triangular decomposition (theorem \ref{dectrian}). Some computations are needed to prove the compatibility with affine quantum Serre relations (section  \ref{proofdecomp}); note that we get a new proof of a combinatorial identity discovered by Jing (consequence of lemma \ref{stepdeux}). The triangular decomposition is used in section \ref{needtrian} to define the Verma modules of $\U_q(\hat{\Glie})$.

\noindent In section \ref{int} we recall the classification of (type $1$) simple integrable highest weight representations of quantum Kac-Moody algebras, and we prove such a classification for quantum affinizations (theorem \ref{simpint}; the proof is analogous to the proof given by Chari-Pressley for quantum affine algebras). The point is to give an adapted definition of a weight which we call a $l$-weight : we need a more precise definition than in the case of quantum affine algebras (a $l$-weight must be characterized by the action of $\U_q(\hat{\Hlie})\subset\U_q(\hat{\Glie})$ on a $l$-weight space). We also give the definition of the category $\mathcal{O}(\U_q(\hat{\Glie}))$. 

\noindent In section \ref{qcar} we construct $q$-characters of integrable modules in the category $\mathcal{O}(\U_q(\hat{\Glie}))$. New technical points are to be considered (in comparison to quantum affine algebra cases) : we have to add terms of the form $k_{\lambda}$ ($\lambda$ coweight of $\U_q(\Glie)$) for the well-definedness in the general case. The original definition of $q$-characters (\cite{Fre}) was based on an explicit formula for the universal $\mathcal{R}$-matrix. In general no universal $\mathcal{R}$-matrix has been defined for a quantum affinization. However $q$-characters can be obtained with a piece of the formula of a ``$\mathcal{R}$-matrix'' in the same spirit as the original approach (theorem \ref{rmat}). In section \ref{comb} we prove that the image of $\chi_q$ is the intersection of the kernels of screening operators (theorem \ref{sym}) in the same spirit as Frenkel-Mukhin \cite{Fre2} did for quantum affine algebras; new technical points are involved because of the $k_{\lambda}$ (we suppose that the quantized Cartan matrix $C(z)$ is invertible). In particular it unifies this approach with \cite{her03} and enables us to prove some results announced in \cite{her03}. We prove that the image of $\chi_q$ is a ring. As $\chi_q$ is injective, we get an induced ring structure $*$ on the Grothendieck group. 

\noindent In section \ref{dnc} we prove that $*$ is a fusion product (theorem \ref{posc}), that is to say that there is a product of modules. We use the new Drinfel'd coproduct (proposition \ref{coprod}); as it involves infinite sums, we have to work in a larger category where the tensor product is well-defined (theorem \ref{prod}). To conclude the proof of theorem \ref{posc} we define specializations of certain forms which allow us to go from the larger category to $\mathcal{O}(\U_q(\hat{\Glie}))$ (section \ref{spe}). We also give some concrete examples of explicit computations in section \ref{ex}.

\noindent {\bf Acknowledgments} : the author would like to thank Marc Rosso for his continued support and Olivier Schiffmann for his accurate remarks.

\section{Background}\label{bck}

\subsection{Cartan matrix} In this section we give some general backgrounds about Cartan matrices (for more details see \cite{kac}). A generalized Cartan matrix is $C=(C_{i,j})_{1\leq i,j\leq n}$ such that\label{carmat} $C_{i,j}\in\ZZ$, $C_{i,i}=2$, $i\neq j\Rightarrow C_{i,j}\leq 0$, $C_{i,j}=0\Leftrightarrow C_{j,i}=0$. We denote $I=\{1,...,n\}$ and $l=\text{rank}(C)$.

\noindent In the following we suppose that $C$ is symmetrizable, that is to say there is a matrix $D=\text{diag}(r_1,...,r_n)$ ($r_i\in\NN^*$)\label{ri} such that $B=DC$\label{symcar} is symmetric. In particular if $C$ is symmetric then it is symmetrizable with $D=I_n$. For example:

\noindent $C$ is said to be of finite type if all its principal minors are in $\NN^*$ (see \cite{bou} for a classification).

\noindent $C$ is said to be of affine type if all its proper principal minor are in $\NN^*$ and $\text{det}(C)=0$ (see \cite{kac} for a classification).

\noindent Let $z$\label{z} be an indeterminate. We put $z_i=z^{r_i}$ and for $l\in\ZZ$, we set $[l]_z=\frac{z^l-z^{-l}}{z-z^{-1}}\in\ZZ[z^{\pm}]$. Let $C(z)$ be the quantized Cartan matrix defined by ($i\neq j\in I$):
$$C_{i,i}(z)=z_i+z_i^{-1}\text{ , }C_{i,j}(z)=[C_{i,j}]_z$$
In sections \ref{comb} and \ref{dnc} we suppose that $C(z)$ is invertible. We have seen in lemma 6.9 of \cite{her03} that the condition $(C_{i,j}<-1\Rightarrow -C_{j,i}\leq r_i)$ implies that $\text{det}(C(z))\neq 0$. In particular finite and affine Cartan matrices (where we impose $r_1=r_2=2$ for $A_1^{(1)}$) satisfy this condition and so the quantum affine algebras and quantum toroidal algebra are included in our study. We denote by $\tilde{C}(z)$ the inverse matrix of $C(z)$ and $D(z)$ the diagonal matrix such that for $i,j\in I$, $D_{i,j}(z)=\delta_{i,j}[r_i]_z$.

\noindent We consider a realization $(\Hlie, \Pi, \Pi^{\vee})$ of $C$ (see \cite{kac}): $\Hlie$ is a $2n-l$ dimensional $\QQ$-vector space, $\Pi=\{\alpha_1,...,\alpha_n\}\subset \Hlie^*$ (set of the simple roots) and $\Pi^{\vee}=\{\alpha_1^{\vee},...,\alpha_n^{\vee}\}\subset \Hlie$ (set of simple coroots) and for $1\leq i,j\leq n$:
$$\alpha_j(\alpha_i^{\vee})=C_{i,j}$$
Denote by $\Lambda_1,...,\Lambda_n\in\Hlie^*$ (resp. the $\Lambda_1^{\vee},...,\Lambda_n^{\vee}\in\Hlie$) the fundamental weights (resp. coweights) : we have $\alpha_i(\Lambda_j^{\vee})=\Lambda_i(\alpha_j^{\vee})=\delta_{i,j}$.

\noindent Consider a symmetric bilinear form $(,):\Hlie^*\times \Hlie^*\rightarrow \QQ$ such that for $i\in I$, $h\in\Hlie^*$ : $(\alpha_i,h)=h(r_i\alpha_i^{\vee})$. It is non degenerate and gives an isomorphism $\nu:\Hlie^*\rightarrow \Hlie$. In particular for $i\in I$ we have $\nu(\alpha_i)=r_i\alpha_i^{\vee}$ and for $\lambda,\mu\in\Hlie^*$, $\lambda(\nu(\mu))=\mu(\nu(\lambda))$.

\noindent Denote $P=\{\lambda \in\Hlie^*/\forall i\in I, \lambda(\alpha_i^{\vee})\in\ZZ\}$ the set of weights and $P^+=\{\lambda \in P/\forall i\in I, \lambda(\alpha_i^{\vee})\geq 0\}$ the set of dominant weights. For example we have $\alpha_1,...,\alpha_n\in P$ and $\Lambda_1,...,\Lambda_n\in P^+$. Denote $Q=\underset{i\in I}{\bigoplus}\ZZ \alpha_i\subset P$ the root lattice and $Q^+=\underset{i\in I}{\sum}\NN \alpha_i\subset Q$. For $\lambda,\mu\in \Hlie^*$, write $\lambda \geq \mu$ if $\lambda-\mu\in Q^+$.

\noindent If $C$ is finite we have $n=l=\text{dim}(\Hlie)$ and for $\lambda\in \Hlie^*$, $\lambda=\underset{i\in I}{\sum}\alpha_i^{\vee}(\lambda)\Lambda_i$. In particular $\alpha_i=\underset{j\in I}{\sum}C_{j,i} \Lambda_j$. In general the simple roots can not be expressed in terms of the fundamental weights.

\subsection{Quantum Kac-Moody algebra}\label{qkma}

\begin{defi} The quantum Kac-Moody algebra $\U_q(\Glie)$ is the $\CC$-algebra with generators $k_h$ ($h\in \Hlie$), $x_i^{\pm}$ ($i\in I$) and relations: 
\begin{equation}\label{cartanqkm}k_hk_{h'}=k_{h+h'}\text{ , }k_0=1\end{equation}
\begin{equation}k_hx_j^{\pm}k_{-h}=q^{\pm \alpha_j(h)}x_j^{\pm}\end{equation}
\begin{equation}[x_i^+,x_j^-]=\delta_{i,j}\frac{k_{r_i\alpha
_i^{\vee}}-k_{-r_i\alpha_i^{\vee}}}{q_i-q_i^{-1}}\end{equation}
\begin{equation}\label{serreqkm}\underset{r=0... 1-C_{i,j}}{\sum}(-1)^r\begin{bmatrix}1-C_{i,j}\\r\end{bmatrix}_{q_i}(x_i^{\pm})^{1-C_{i,j}-r}x_j^{\pm}(x_i^{\pm})^r=0 \text{ (for $i\neq j$)}\end{equation}
\end{defi}

\noindent This algebra was introduced independently by Jimbo \cite{jim} and Drinfel'd \cite{Dri1} and is also called a quantum group. It is remarkable that one can define a Hopf algebra structure on $\U_q(\Glie)$ by setting :
$$\Delta(k_h)=k_h\otimes k_h$$
$$\Delta(x_i^+)=x_i^+\otimes 1 + k_i^+\otimes x_i^+\text{ , }\Delta(x_i^-)=x_i^-\otimes k_i^- + 1\otimes x_i^-$$
$$S(k_h)=k_{-h}\text{ , }S(x_i^+)=-x_i^+k_i^{-1}\text{ , }S(x_i^-)=-k_i^+x_i^-$$
$$\epsilon(k_h)=1\text{ , }\epsilon(x_i^+)=\epsilon(x_i^-)=0$$
where we use the notation $k_i^{\pm}=k_{\pm r_i\alpha_i^{\vee}}$. 

\noindent For $i\in I$ let $U_i$ be the subalgebra of $\U_q(\Glie)$ generated by the $x_i^{\pm}, k_{p\alpha_i^{\vee}}$ ($p\in\QQ$). Then $U_i$ is isomorphic to $\U_{q_i}(sl_2)$, and so a $\U_q(\Glie)$-module has also a structure of $\U_{q_i}(sl_2)$-module.

\begin{defi} A triangular decomposition of an algebra $A$ is the data of three subalgebras $(A^-, H, A^+)$ of $A$ such that the multiplication  $x^-\otimes h\otimes x^+\mapsto x^-hx^+$ defines an isomorphism of $\CC$-vector space $A^-\otimes H\otimes A^+\simeq A$.\end{defi}

\noindent Let $\U_q(\Glie)^+$ (resp. $\U_q(\Glie)^-$, $\U_q(\Hlie)$) be the subalgebra of $\U_q(\Glie)$ generated by the $x_i^+$ (resp. the $x_i^-$, resp. the $k_h$). We have (see \cite{lu}) : 

\begin{thm}\label{triansc} $(\U_q(\Glie)^-,\U_q(\Hlie),\U_q(\Glie)^+)$ is a triangular decomposition of $\U_q(\Glie)$. Moreover $\U_q(\Hlie)$ (resp. $\U_q(\Glie)^+$, $\U_q(\Glie)^-$) is isomorphic to the algebra with generators $k_h$ (resp $x_i^+$, $x_i^-$) and relations (\ref{cartanqkm}) (resp. relations (\ref{serreqkm}) with $+$, relations (\ref{serreqkm}) with $-$).\end{thm}

\section{Quantum affinization $\U_q(\hat{\Glie})$ and triangular decomposition}\label{genaff}

\noindent In this section we define general quantum affinizations (without central charge), we give the relations between the currents (section \ref{cur}) and we prove a triangular decomposition (theorem \ref{dectrian}).

\subsection{Definition}

\begin{defi}\label{defiaffi} The quantum affinization of $\U_q(\Glie)$ is the $\CC$-algebra $\U_q(\hat{\Glie})$ with generators $x_{i,r}^{\pm}$ ($i\in I, r\in\ZZ$), $k_h$ ($h\in \Hlie$), $h_{i,m}$ ($i\in I, m\in\ZZ-\{0\}$) and the following relations ($i,j\in I, r,r'\in\ZZ, m\in\ZZ-\{0\}$):
\begin{equation}\label{afcart}k_hk_{h'}=k_{h+h'}\text{ , }k_0=1\text{ , }[k_{h},h_{j,m}]=0\text{ , }[h_{i,m},h_{j,m'}]=0\end{equation}
\begin{equation}\label{actcartlu}k_{h}x_{j,r}^{\pm}k_{-h}=q^{\pm \alpha_j(h)}x_{j,r}^{\pm}\end{equation}
\begin{equation}\label{actcartld}[h_{i,m},x_{j,r}^{\pm}]=\pm \frac{1}{m}[mB_{i,j}]_qx_{j,m+r}^{\pm}\end{equation}
\begin{equation}\label{pml}[x_{i,r}^+,x_{j,r'}^-]= \delta_{ij}\frac{\phi^+_{i,r+r'}-\phi^-_{i,r+r'}}{q_i-q_i^{-1}}\end{equation}
\begin{equation}\label{equaun}x_{i,r+1}^{\pm}x_{j,r'}^{\pm}-q^{\pm B_{ij}}x_{j,r'}^{\pm}x_{i,r+1}^{\pm}=q^{\pm B_{ij}}x_{i,r}^{\pm}x_{j,r'+1}^{\pm}-x_{j,r'+1}^{\pm}x_{i,r}^{\pm}\end{equation}
\begin{equation}\label{equadeux}\underset{\pi\in \Sigma_s}{\sum}\underset{k=0..s}{\sum}(-1)^k\begin{bmatrix}s\\k\end{bmatrix}_{q_i}x_{i,r_{\pi(1)}}^{\pm}...x_{i,r_{\pi(k)}}^{\pm}x_{j,r'}^{\pm}x_{i,r_{\pi(k+1)}}^{\pm}...x_{i,r_{\pi(s)}}^{\pm}=0\end{equation}
where the last relation holds for all $i\neq j$, $s=1-C_{ij}$, all sequences of integers $r_1,...,r_s$. $\Sigma_s$ is the symmetric group on $s$ letters. For $i\in I$ and $m\in\ZZ$, $\phi_{i,m}^{\pm}\in \U_q(\hat{\Glie})$ is determined by the formal power series in $\U_q(\hat{\Glie})[[z]]$ (resp. in $\U_q(\hat{\Glie})[[z^{-1}]]$):
$$\underset{m\geq 0}{\sum}\phi_{i,\pm m}^{\pm}z^{\pm m}=k_{\pm r_i \alpha_i^{\vee}}\text{exp}(\pm(q-q^{-1})\underset{m'\geq 1}{\sum}h_{i,\pm m'}z^{\pm m'})$$
and $\phi_{i,m}^+=0$ for $m<0$, $\phi_{i,m}^-=0$ for $m>0$.
\end{defi}

\noindent The relations (\ref{equadeux}) are called affine quantum Serre relations. The notation $k_i^{\pm}=k_{\pm r_i\alpha_i^{\vee}}$ is also used. We have $k_ik_i^{-1}=k_i^{-1}k_i=1\text{ , }k_ix_{j,m}^{\pm}k_i^{-1}=q^{\pm B_{ij}}x_{j,m}^{\pm}$.

\noindent There is an algebra morphism $\U_q(\Glie)\rightarrow \U_q(\hat{\Glie})$ defined by ($h\in \Hlie, i\in I$) $k_h\mapsto k_h$ , $x_i^{\pm}\mapsto x_{i,0}^{\pm}$. In particular a $\U_q(\hat{\Glie})$-module has also a structure of a $\U_q(\Glie)$-module.

\subsection{Relations between the currents}\label{cur} For $i\in I$, consider the series (also called currents):
$$x_i^{\pm}(w)=\underset{r\in\ZZ}{\sum}x_{i,r}^{\pm}w^r\text{ , }\phi_i^+(z)=\underset{m\geq 0}{\sum}\phi_{i,m}^+z^m\text{ , }\phi_i^-(z)=\underset{m\geq 0}{\sum}\phi_{i,-m}^-z^{-m}$$
The defining relations of $\U_q(\hat{\Glie})$ can be written with currents ($h,h'\in\Hlie$, $i,j\in I$):
\begin{equation}\label{afcartu}k_hk_{h'}=k_{h+h'}\text{ , }k_0=1\text{ , }k_h\phi_i^{\pm}(z)=\phi_i^{\pm}(z)k_h\end{equation}
\begin{equation}\label{actcartc}k_hx_j^{\pm}(z)=q^{\pm \alpha_j (h)}x_j^{\pm}(z)k_h\end{equation}
\begin{equation}\label{actcartplus}\phi_i^+(z)x_j^{\pm}(w)=\frac{q^{\pm B_{i,j}}w-z}{w-q^{\pm B_{i,j}}z}x_j^{\pm}(w)\phi_i^+(z)\end{equation}
\begin{equation}\label{actcartmoins}\phi_i^-(z)x_j^{\pm}(w)=\frac{q^{\pm B_{i,j}}w-z}{w-q^{\pm B_{i,j}}z}x_j^{\pm}(w)\phi_i^-(z)\end{equation}
\begin{equation}\label{partact}[x_i^+(z),x_j^-(w)]=\frac{\delta_{i,j}}{q_i-q_i^{-1}}[\delta(\frac{w}{z})\phi_i^+(w)-\delta(\frac{z}{w})\phi_i^-(z)]\end{equation}
\begin{equation}\label{plusmoinsc}(w-q^{\pm B_{i,j}}z)x_i^{\pm}(z)x_j^{\pm}(w)=(q^{\pm B_{i,j}}w-z)x_j^{\pm}(w)x_i^{\pm}(z)\end{equation}
\begin{equation}\label{equadeuxc}\underset{\pi\in \Sigma_s}{\sum}\underset{k=0...s}{\sum}(-1)^k\begin{bmatrix}s\\k\end{bmatrix}_{q_i}x_i^{\pm}(w_{\pi(1)})...x_i^{\pm}(w_{\pi(k)})x_j^{\pm}(z)x_i^{\pm}(w_{\pi(k+1)})...x_i^{\pm}(w_{\pi(s)})=0\end{equation}
where $\delta(z)=\underset{r\in\ZZ}{\sum}z^r$. The equation (\ref{actcartplus}) (resp. equation (\ref{actcartmoins})) is expanded for $|z|<|w|$ (resp. $|w|<|z|$). 

\noindent Remark: in the relations (\ref{plusmoinsc}), the terms can not be divided by $w-q^{\pm B_{i,j}}z$ : it would involve infinite sums and make no sense.

\noindent The following equivalences are clear : (relations (\ref{afcart}) $\Leftrightarrow$ relations (\ref{afcartu})) ; (relations (\ref{actcartlu}) $\Leftrightarrow$ relations (\ref{actcartc})) ; (relations (\ref{equaun}) $\Leftrightarrow$ relations (\ref{plusmoinsc})) ; (relations (\ref{pml}) $\Leftrightarrow$ relations (\ref{partact})) ; (relations (\ref{equadeux}) $\Leftrightarrow$ relations (\ref{equadeuxc})).

\noindent We suppose that the relations (\ref{actcartlu}) are verified and we prove the equivalence (relations (\ref{actcartld}) with $m\geq 1$ $\Leftrightarrow$ relations (\ref{actcartplus})) ((relations (\ref{actcartld}) with $m\leq -1$ $\Leftrightarrow$ relations (\ref{actcartmoins})) is proved in an similar way): consider $h_i^+(z)=\underset{m\geq 1}{\sum}mh_{i,m}z^{m-1}$. The relation (\ref{actcartld}) with $m\geq 1$ are equivalent to (expanded for $|z|<|w|$):
$$[h_i^+(z),x_j^{\pm}(w)]=\pm [B_{i,j}]_q\frac{w^{-1} x_j^{\pm}(w)}{(1-\frac{z}{w}q^{B_{i,j}})(1-\frac{z}{w}q^{-B_{i,j}})}$$
It is equivalent to the data of a $\alpha_{\pm}(z,w)\in(\CC[w, w^{-1}])[[z]]$ such that $\phi_i^+(z)x_j^{\pm}(w)=\alpha_{\pm}(z,w)x_j^{\pm}(w)\phi_i^+(z)$. So it suffices to prove that this term is the $\frac{q^{\pm B_{i,j}}w-z}{w-q^{\pm B_{i,j}}z}$ of relation (\ref{actcartplus}). Let us compute this term : we have $\frac{\partial \phi_i^+(z)}{\partial z}=(q-q^{-1})h_i^+(z)\phi_i^+(z)$ and so the relations (\ref{actcartld}) imply :
$$(q-q^{-1})\phi_i^+(z)[h_i^+(z),x_j^{\pm}(w)]=\frac{\partial \alpha_{\pm}(z,w)}{\partial z}x_j^{\pm}(w)\phi_i^+(z)$$
$$(\pm [B_{i,j}]_q\frac{w^{-1}}{(1-\frac{z}{w}q^{B_{i,j}})(1-\frac{z}{w}q^{-B_{i,j}})}\alpha_{\pm}(z,w)-\frac{1}{q-q^{-1}}\frac{\partial \alpha_{\pm}(z,w)}{\partial z})x_j^{\pm}(w)\phi_i^+(z)=0$$
$$\frac{\partial \alpha_{\pm}(z,w)}{\partial z}=\pm (q^{B_{i,j}}-q^{-B_{i,j}})\frac{w^{-1}}{(1-\frac{z}{w}q^{B_{i,j}})(1-\frac{z}{w}q^{-B_{i,j}})}\alpha_{\pm}(z,w)$$
As $\frac{q^{\pm B_{i,j}}w-z}{w-q^{\pm B_{i,j}}z}$ is a solution, we have $\alpha_{\pm}(z,w)=\lambda(w)\frac{q^{\pm B_{i,j}}w-z}{w-q^{\pm B_{i,j}}z}$. But at $z=0$ we know $\alpha_{\pm}(0,w)=q^{\pm B_{i,j}}$ (relations (\ref{actcartlu})) and so $\lambda(w)=1$.

\subsection{Triangular decomposition}\label{proofdecomp}

\subsubsection{Statement}\label{stat} Let $\U_q(\hat{\Glie})^+$ (resp. $\U_q(\hat{\Glie})^-$, $\U_q(\hat{\Hlie})$) be the subalgebra of $\U_q(\hat{\Glie})$ generated by the $x_{i,r}^+$ (resp. the $x_{i,r}^-$, resp. the $k_h$, $h_{i,r}$).

\begin{thm}\label{dectrian} $(\U_q(\hat{\Glie})^-, \U_q(\hat{\Hlie}), \U_q(\hat{\Glie})^+)$ is a triangular decomposition of $\U_q(\hat{\Glie})$. Moreover $\U_q(\hat{\Hlie})$ (resp. $\U_q(\hat{\Glie})^+$, $\U_q(\hat{\Glie})^-$) is isomorphic to the algebra with generators $k_h, h_{i,m}$ (resp $x_{i,r}^+$, $x_{i,r}^-$) and relations (\ref{afcart}) (resp. relations (\ref{equaun}), (\ref{equadeux}) with $+$, relations (\ref{equaun}), (\ref{equadeux}) with $-$).\end{thm}

\noindent  For a quantum affine algebra ($C$ finite) it is proved in \cite{bec}. 

\noindent In this section \ref{proofdecomp} we prove this theorem in general. We will use the algebras $\U_q^l(\hat{\Glie})$, $\tilde{\U}_q(\hat{\Glie})$ defined by :

\begin{defi} $\U_q^l(\hat{\Glie})$ is the $\CC$-algebra with generators $x_{i,r}^{\pm}$, $h_{i,m}$, $k_h$ ($i\in I$, $r\in\ZZ$, $m\in\ZZ-\{0\}$, $h\in\Hlie$) and relations (\ref{afcart}), (\ref{actcartlu}), (\ref{actcartld}), (\ref{pml}) (or relations (\ref{afcartu}), (\ref{actcartc}), (\ref{actcartplus}), (\ref{actcartmoins}), (\ref{partact})). 

\noindent $\tilde{\U}_q(\hat{\Glie})$ is the quotient of $\U_q^l(\hat{\Glie})$ by relations (\ref{equaun}) (or relations (\ref{plusmoinsc})).\end{defi}

\noindent Note that $\U_q(\hat{\Glie})$ is a quotient of $\U_q^l(\hat{\Glie})$ and that $(\U_q^{l,-}(\hat{\Glie}), \U_q(\hat{\Hlie}), \U_q^{l,+}(\hat{\Glie}))$ is a triangular decomposition of $\U_q^l(\hat{\Glie})$ where $\U_q^{l,\pm}(\hat{\Glie})$ is generated by the $x_{i,r}^{\pm}$ without relations. In the $sl_2$-case we have $\tilde{\U}_q(\hat{sl_2})=\U_q(\hat{sl_2})$. 

\noindent Let us sketch the proof of theorem \ref{dectrian}. We use a method analog to the proof for classic cases or quantum Kac-Moody algebras (see for example the chapter 4 of \cite{jan}) : we have to check a compatibility condition between the relations and the product as explained in section \ref{humgen}. After some preliminary technical lemmas about polynomials in section \ref{plustech}, the heart of the proof is given in section \ref{finpreuve} : properties of $\U_q^l(\hat{\Glie})$ (lemma \ref{stepun}) lead to a triangular decomposition of $\tilde{\U}_q(\hat{\Glie})$. Properties of $\tilde{\U}_q(\hat{\Glie})$ proved in lemmas \ref{stepdeux}, \ref{steptrois} imply theorem \ref{dectrian}. Note that the intermediate algebra $\tilde{\U}_q(\hat{\Glie})$ is also studied because it will be used in the last section of this paper.

\noindent Remark : lemma \ref{stepdeux} gives a new proof of a combinatorial identity discovered by Jing.

\noindent The theorem \ref{dectrian} is used in section \ref{needtrian} to define the Verma modules of $\U_q(\hat{\Glie})$. Let us give another consequence of theorem \ref{dectrian} : for $i\in I$, let $\hat{U}_i$ be the subalgebra of $\U_q(\hat{\Glie})$ generated by the $x_{i,r}^{\pm}, k_{p\alpha_i^{\vee}}, h_{i,m}$ ($r\in\ZZ$, $m\in\ZZ-\{0\}$, $p\in\QQ$). We have a morphism $\U_{q_i}(\hat{sl_2})\rightarrow \hat{U}_i$ (in particular any $\U_q(\hat{\Glie})$-module also has a structure of $\U_{q_i}(\hat{sl_2})$-module). Moreover theorem \ref{dectrian} implies:

\begin{cor} $\hat{U}_i$ is isomorphic to $\U_{q_i}(\hat{sl_2})$.\end{cor}

\subsubsection{General proof of triangular decompositions}\label{humgen} Let $A$ be an algebra with a triangular decomposition $(A^-,H,A^+)$. Let $B^+$ (resp. $B^-$) be a two-sided ideal of $A^+$ (resp. $A^-$). Let $C=A/(A.(B^+ + B^-).A)$ and denote by $C^{\pm}$ the image of $B^{\pm}$ in $C$.

\begin{lem}\label{gendec} If $B^+.A\subset A.B^+$ and $A.B^-\subset B^-.A$ then $(C^-,H,C^+)$ is a triangular decomposition of $C$ and the algebra $C^{\pm}$ is isomorphic to $A^{\pm}/B^{\pm}$.\end{lem}

\demo We use the proof of section 4.21 in \cite{jan} : indeed the product gives an isomorphism of $\CC$-vector space $A.(B^+ + B^-).A\simeq B^+\otimes H\otimes A^- + A^+\otimes H\otimes B^-$.\qed

\subsubsection{Technical lemmas}\label{plustech} Let $i\neq j$ and $s=1-C_{i,j}$. Define $P_{\pm}(w_1,...,w_s,z)\in\CC[w_1,...,w_s,z]$ by the formula :
$$P_{\pm}(w_1,...,w_s,z)=\underset{k=0...s}{\sum}(-1)^k\begin{bmatrix}s\\k\end{bmatrix}_{q_i}(w_1-q^{\pm B_{i,j}}z)...(w_k-q^{\pm B_{i,j}}z)(w_{k+1}q^{\pm B_{i,j}}-z)...(w_sq^{\pm B_{i,j}}-z)$$

\begin{lem}\label{techun} There are polynomials $(f_{\pm,r})_{r=1,...,s-1}$ of $s-1$ variables such that:
$$P_{\pm}(w_1,...,w_s,z)=\underset{1\leq r\leq s-1}{\sum}(w_{r+1}-q_i^{\pm 2} w_r)f_{\pm,r}(w_1,...,w_{r-1},w_{r+2},...,w_s,z)$$\end{lem}

\demo It suffices to prove it for $P_+$ (because $P_-$ is obtained from $P_+$ by $q\mapsto q^{-1}$). First we prove that $P_+(q_i^{-2(s-1)}w,q_i^{-2(s-2)}w,...,q_i^{-2}w,w,z)=0$. Indeed it is equal to: 
$$w^s\underset{k=0...s}{\sum}(-1)^k\begin{bmatrix}s\\k\end{bmatrix}_{q_i}q_i^{k (1-s)}(q_i^{-2(s-1)}q_i^{s-1}-\frac{z}{w})...(q_i^{-2(s-k)}q_i^{s-1}-\frac{z}{w})(q_i^{-2(s-k-1)}q_i^{1-s}-\frac{z}{w})...(q_i^{1-s}-\frac{z}{w})$$
$$=z^s(q_i^{1-s}-\frac{z}{w})q^{-3s+3}M_{q_i}(\frac{z}{w}q^{3s-3})$$
where:
$$M_q(u)=\underset{k=0...s}{\sum}(-1)^k\begin{bmatrix}s\\k\end{bmatrix}_qq^{k(1-s)}(q^{2k}-u)(q^{2(k+1)}-u)...(q^{2(k+s-2)}-u)$$
Let $\alpha_0(q),...,\alpha_{s-1}(q)\in\ZZ[q]$ such that $(a-u)(a-uq^2)...(a-uq^{2(s-2)})=u^{s-1}\alpha_{s-1}(q)+u^{s-2}a\alpha_{s-2}(q)+...+a^{s-1}\alpha_0(q)$. So:
$$M_q(u)=\underset{p=0...s-1}{\sum}\alpha_{s-p}(q)u^{s-p}\underset{k=0...s}{\sum}(-1)^k\begin{bmatrix}s\\k\end{bmatrix}_qq^{k(1-s+2p)}$$
And so $M_q(u)=0$ because of the $q$-binomial identity for $p'=1-s,3-s,...,s-1$ (see \cite{lu}):
$$\underset{k=0...s}{\sum}(-1)^k\begin{bmatrix}s\\k\end{bmatrix}_qq^{rp'}=0$$
As a consequence $P_+$ is in the kernel of the projection 
$$\phi : \CC[w_1,...,w_s,z]\rightarrow \CC[w_1,...,w_s,z]/((w_2-q_i^2 w_1),...,(w_s-q_i^2w_{s-1}))$$ 
that is to say $P_+(w_1,...,w_s,z)=\underset{1\leq r\leq s-1}{\sum}(w_{r+1}-q^{B_{i,j}}w_r)f_r(w_1,...,w_s,z)$ where the $f_r\in\CC[w_1,...,w_s,z]$. 

\noindent Let us prove that we can choose the $(f_r)_{1\leq r\leq s-1}$ so that for all $1\leq s\leq r-1$, $f_r$ does not depend of $w_r, w_{r+1}$. Let $\mathcal{A}\subset \text{Ker}(\phi)$ be the subspace of polynomials which are degree at most of $1$ in each variable $w_1,...,w_s$. In particular $P\in\mathcal{A}$. We can decompose in a unique way $P=\alpha+w_2\beta+w_1\gamma$ where $\alpha,\gamma\in\CC[w_3,...,w_s,z]$, $\beta\in\CC[w_1,w_3,...,w_s,z]$. Consider $\lambda^{(1)}=-q_i^{-2}\gamma(w_2-q_i^2w_1)\in\mathcal{A}$ and $P^{(1)}=P-\lambda^{(1)}\in\mathcal{A}$. We have in particular $P^{(1)}=\mu^{(1)}_3 + w_2 \mu^{(1)}_2 + w_2 w_1 \mu^{(1)}_1$ where $\mu^{(1)}_1, \mu^{(1)}_2, \mu^{(1)}_3 \in\CC[w_3,...,w_s,z]$. In the same way we define by induction on $r$ ($1\leq r\leq s-1$) the $\lambda^{(r)}\in\mathcal{A}$ such that $P^{(r)}=P^{(r-1)}-\lambda^{(r)}\in \mathcal{A}$ is of the form:
$$P^{(r)}=\mu^{(r)}_{r+2} + w_{r+1}\mu^{(r)}_{r+1} + w_{r+1}w_r \mu^{(r)}_r+...+ w_{r+1}w_r...w_1 \mu^{(r)}_1$$
where for $1\leq r'\leq r+2$, $\mu^{(r)}_{r'}\in\CC[w_{r+2},...,w_s,z]$. Indeed in the part of $P^{(r)}$ without $w_{r+2}$ we can change the terms $w_{r+1}\lambda(w_{r+3},...,w_s,z)$ to $q_i^{-2}w_{r+2}\lambda(w_{r+3},...,w_s,z)$ by adding $q_i^{-2}(w_{r+2}-q_i^2w_{r+1})\lambda\in\mathcal{A}$, we can change the terms $w_{r+1}w_r\lambda'(w_{r+3},...,w_s,z)$ to $q_i^{-4}w_{r+2}w_{r+1}\lambda'(w_{r+3},...,w_s,z)$ by adding $q_i^{-4}(w_{r+2}-q_i^2w_{r+1})\lambda+q_i^{-2}(w_{r+2}-q_i^2w_{r+1})\lambda\in\mathcal{A}$, and so on. In particular for $r=s-1$ :
$$P^{(s-1)}=\mu^{(s-1)}_{s+1} +\mu^{(s-1)}_s w_s+\mu^{(s-1)}_{s-1}w_sw_{s-1}+...+\mu^{(s-1)}_1 w_sw_{s-1}...w_1$$
where $\mu^{(s-1)}_{s+1},...,\mu^{(s-1)}_1\in\CC[z]$. But :
$$0=\phi(P^{(s-1)})=\mu^{(s-1)}_{s+1} +\mu^{(s-1)}_s w_s+\mu^{(s-1)}_{s-1}q_i^{-2}w_s^2+...+\mu^{(s-1)}_1 q_i^{-2-4-...-2(s-1)}w_s^s$$
So for all $1\leq r'\leq s+1$, $\mu^{(s-1)}_{r'}=0$, and so $P^{(s-1)}=0$. In particular $P=\lambda^{(1)}+\lambda^{(2)}+...+\lambda^{(s-1)}$.\qed

For $1\leq k \leq s$ consider $P_{\pm}^{(k)}(w_1,w_2,...,w_s,z)\in\CC[w_1,...,w_s,z]$ defined by:
$$(-1)^k\begin{bmatrix}s\\k\end{bmatrix}_{q_i}\underset{k'=1...k}{\sum}(zq_i^{\pm (1-s)}-w_1)(w_2-q_i^{\pm 2} w_1)...(w_{k'}-q_i^{\pm 2}w_1)(w_{k'+1}q_i^{\pm 2}-w_1)...(w_s q_i^{\pm 2}-w_1)$$
$$+(-1)^{k-1}\begin{bmatrix}s\\k-1\end{bmatrix}_{q_i}\underset{k'=k...s}{\sum}(z-w_1q_i^{\pm (1-s)})(w_2-q_i^{\pm 2} w_1)...(w_{k'}-q_i^{\pm 2}w_1)(w_{k'+1}q_i^{\pm 2}-w_1)...(w_sq_i^{\pm 2}-w_1)$$

\begin{lem}\label{techdeux} i) For $2\leq k\leq s-1$ there are polynomials $(f_{\pm,r}^{(k)})_{r=1,...,s-1}$ of $s-1$ variables, of degree at most $1$ in each variable, such that $P_{\pm}^{(k)}(w_1,...,w_s,z)$ is equal to :
$$(z-q_i^{\pm (1-s)}w_k)f_{\pm ,k-1}^{(k)}(w_1,...,w_{k-1},w_{k+1},...,w_s,z)+(w_{k+1}-q_i^{\pm (1-s)}z)f_{\pm, s-1}^{(k)}(w_1,...,w_k,w_{k+2},...,w_s,z)$$
$$+\underset{1\leq r\leq s-2, r\neq k-1}{\sum}(w_{r+2}-q_i^{\pm 2} w_{r+1})f_{\pm,r}^{(k)}(w_1,...,w_{r-1},w_{r+2},...,w_s,z)$$

ii) There are polynomials $(f_{\pm,r}^{(1)})_{r=1,...,s-1}$ of $s-1$ variables, of degree at most $1$ in each variable, such that $P_{\pm}^{(1)}(w_1,...,w_s,z)$ is equal to :
$$(w_2-q_i^{\pm (1-s)}z)f_{\pm, s}^{(1)}(w_3,...,w_s,z)+\underset{1\leq r\leq s-2}{\sum}(w_{r+2}-q_i^{\pm 2} w_{r+1})f_{\pm,r}^{(k)}(w_1,...,w_{r-1},w_{r+2},...,w_s,z)$$

iii) There are polynomials $(f_{\pm,r}^{(s)})_{r=1,...,s-1}$ of $s-1$ variables, of degree at most $1$ in each variable, such that $P_{\pm}^{(s)}(w_1,...,w_s,z)$ is equal to :
$$(z-q_i^{\pm (1-s)}w_s)f_{\pm, s-1}^{(s)}(w_1,...,w_{s-1},z)+\underset{1\leq r\leq s-2}{\sum}(w_{r+2}-q_i^{\pm 2} w_{r+1})f_{\pm,r}^{(s)}(w_1,...,w_{r-1},w_{r+2},...,w_s,z)$$\end{lem}

\demo It suffices to prove it for $P_+^{(k)}$ (because $P_-^{(k)}$ is obtained from $P_+^{(k)}$ by $q\mapsto q^{-1}$). 

\noindent For i) : we see as in lemma \ref{techun} that it suffices to check that $P_+^{(k)}(w_1,...,w_s,z)=0$ if $w_3=q_i^2w_2$, ... , $w_k=q_i^2w_{k-1}$, $w_{k+2}=q_i^2w_{k+1}$, ... , $w_s=q_i^2w_{s-1}$, $z=q_i^{1-s}w_k$ and $w_{k+1}=q_i^{1-s}z$. It means $w_3=q_i^2w_2$, ..., $w_k=q_i^{2(k-2)}w_2$, $w_{k+1}=q^{2k-2-2s}w_2$, ... ,$w_s=q^{-4}w_2$, $z=q^{2k-3-s} w_2$. So if we set $u=w_1/w_2$ we find for $P_+^{(k)}w_2^{-s}$:

\noindent $(-1)^k\begin{bmatrix}s\\k\end{bmatrix}_{q_i}\underset{k'=1...k}{\sum}q_i^{2(k'-1)}(q_i^{2k-2-2s}-u)(q_i^{2k-2s}-u)...(q_i^{2k'-6}-u)(q_i^{2k'}-u)...(q_i^{2k-2}-u)(q_i^{-2}-u)
\\+(-1)^{k-1}\begin{bmatrix}s\\k-1\end{bmatrix}_{q_i}\underset{k'=k...s}{\sum}q_i^{2k'-s-1}(q_i^{2k-2s-4}-u)...(q_i^{2k'-2s-6}-u)(q_i^{2k'-2s}-u)...(q_i^{2k-4}-u)(q_i^{-2}-u)$

\noindent It is a multiple of:

\noindent $\frac{[s-k+1]_{q_i}}{q_i^{2k-2s-4}-u}[\underset{k'=1...k}{\sum}\frac{q_i^{2k'-1}}{(q_i^{2k'-2}-u)(q_i^{2k'-4}-u)}]-\frac{[k]_{q_i}}{q_i^{2k-2}-u}q_i^s[\underset{k'=k...s}{\sum}\frac{q_i^{2k'-1-s}}{(q_i^{2k'-2s-2}-u)(q_i^{2k'-2s-4}-u)}]
\\=\frac{q_i^2[s-k+1]_{q_i}}{(1-q_i^2)(q_i^{2k-2s-4}-u)}[\underset{k'=1...k}{\sum}\frac{1}{q_i^{2k'-2}-u}-\frac{1}{q_i^{2k'-4}-u}]-\frac{q_i^2 [k]_{q_i}}{(1-q_i^2)q_i^{2k-2}-u}q_i^s[\underset{k'=k...s}{\sum}\frac{1}{q_i^{2k'-2s-2}-u}-\frac{1}{q_i^{2k'-2s-4}-u}]
\\=\frac{q_i^2[s-k+1]_{q_i}}{(1-q_i^2)(q_i^{2k-2s-4}-u)}[\frac{1}{q_i^{2k-2}-u}-\frac{1}{q_i^{-2}-u}]-\frac{q_i^2[k]_{q_i}}{(1-q_i^2)(q_i^{2k-2}-u)}q_i^s[\frac{1}{q_i^{-2}-u}-\frac{1}{q_i^{2k-2s-4}-u}]=0$

\noindent For ii) : as for i) we check that $P_+^{(1)}(w_1,...,w_s,z)=0$ if $w_3=q_i^2w_2$, ... , $w_s=q_i^2w_{s-1}$, $z=q_i^{s-1}w_2$. It means $w_{k'}=q_i^{2(k'-2)}w_2$ for $2\leq k'\leq s$. So it we set $u=w_1/w_2$ we find for $P_+^{(1)}w_2^{-s}$:

\noindent $-[s]_{q_i}(1-u)(q_i^2-u)...(q_i^{2s-2}-u)
+q_i^{1-s}\underset{k'=1...s}{\sum}q_i^{2k'-2}(q_i^{2s-2}-u)(q_i^{-2}-u)...(q_i^{2k'-6}-u)(q_i^{2k'}-u)...(q_i^{2s-2}-u)$

\noindent It is a multiple of: $-\frac{q_i^2[s]_{q_i}}{q_i^{2s}-1}(\frac{1}{q_i^{-2}-u}-\frac{1}{q_i^{2s-2}-u})+\frac{q_i^{1-s}}{1-q_i^{-2}}(\frac{1}{q_i^{-2}-u}-\frac{1}{q_i^{2s-2}-u})=0$.

\noindent For iii) : as for i) we check that $P_+^{(k)}(w_1,...,w_s,z)=0$ if $w_3=q_i^2w_2$, ... , $w_s=q_i^2w_{s-1}$, $z=q_i^{1-s}w_s$. It means $w_{k'}=q_i^{2(k'-2)}w_2$ for $2\leq k'\leq s$ and $z=q_i^{s-3}w_2$. The computation is analogous to i).\qed

\begin{lem}\label{techtrois} For all choices of polynomials $(f_{\pm,r}^{(k')})_{1\leq k'\leq s, 1\leq r\leq s-1}$ in lemma \ref{techdeux} and each $2\leq k\leq s$ there are polynomials $(g_{\pm,r}^{(k)})_{r=1,...,s-2}$ of $s-1$ variables such that:
$$f_{\pm , k-1}^{(k)}-f_{\pm , s-1}^{(k-1)}=\underset{1\leq r\leq s-2}{\sum}(w_{r+2}-q_i^{\pm 2} w_{r+1})g_{\pm,r}^{(k)}(w_1,...,w_{r-1},w_{r+2},...,w_s,z)$$\end{lem}

\demo We see as in lemma \ref{techun} that it suffices to check that $f_{+,k-1}^{(k)}+f_{+,s-1}^{(k-1)}=0$ if $w_3=q_i^2w_2$, ... , $w_s=q_i^2w_{s-1}$. So we suppose that $w_{k'}=q_i^{2(k'-2)}$ for all $2\leq k'\leq s$. Let $Q=w_1^{s-1}(w_2q_i^{-2}-1)(w_2-1)...(w_2q_i^{2s-2})/(q_i^2-1)$). It suffices to prove that for $2\leq k\leq s$, we have:
\begin{equation}\label{trouvun}(q_i^2-q_i^{2-2s})f_{+,k-1}^{(k)}(Q(-1)^k\begin{bmatrix}s\\k\end{bmatrix}_{q_i}[k]_{q_i}(q_i-q_i^{-1}))^{-1}=\frac{q_i^{k+1+s}+q_i^{-s-k+3}-q_i^{-s+k+1}-q_i^{3-k+s}}{(vq_i^{-2}-1)(vq_i^{2k-4}-1)(vq_i^{2s-2}-1)}\end{equation}
\begin{equation}\label{trouvdeux}(q_i^2-q_i^{2-2s})f_{+,s-1}^{(k-1)}(Q(-1)^{k-1}\begin{bmatrix}s\\k-1\end{bmatrix}_{q_i}[k-1]_{q_i}(q_i-q_i^{-1}))^{-1}=\frac{q_i^{k+1}+q_i^{-k+3}-q_i^{-2s+k+1}-q_i^{3-k+2s}}{(vq_i^{-2}-1)(vq_i^{2k-4}-1)(vq_i^{2s-2}-1)}\end{equation}
because we have the relation :
$$\begin{bmatrix}s\\k\end{bmatrix}_{q_i}[k]_{q_i}(q_i^{k+1+s}+q_i^{-s-k+3}-q_i^{-s+k+1}-q_i^{3-k+s})=-\begin{bmatrix}s\\k-1\end{bmatrix}_{q_i}[k-1]_{q_i}(q_i^{k+1}+q_i^{-k+3}-q_i^{-2s+k+1}-q_i^{3-k+2s})$$
First suppose that $3\leq k\leq s-1$. We have $P_+^{(k)}=(z-q_i^{1-s}w_k)f_{+,k-1}^{(k)}+(q_i^2w_k-q_i^{1-s}z)f_{+,s-1}^{(k)}$. So for $\alpha_k,\beta_k$ such that $P_+^{(k)}=z\alpha_k+w_k \beta_k$, we have $f_{+,k-1}^{(k)}=\frac{q_i^2\alpha_k+q_i^{1-s}\beta_k}{q_i^2-q_i^{2-2s}}$ and $f_{+,s-1}^{(k)}=\frac{q_i^{1-s}\alpha_k+\beta_k}{q_i^2-q_i^{2-2s}}$. But we have $P_+^{(k)}=z(q_i^{1-s}\lambda_k+\mu_k)-w_1(\lambda_k+q_i^{1-s}\mu_k)$ where (we put $v=w_2/w_1$):

\noindent $\lambda_k=(-1)^kw_1^{s-1}\begin{bmatrix}s\\k\end{bmatrix}_{q_i}\underset{k'=1...k}{\sum}(v-q_i^2)(vq_i^2-q_i^2)...(vq_i^{2(k'-2)}-q_i^2)(vq_i^{2k'}-1)...(vq_i^{2(s-2)+2}-1)
\\=Q(-1)^k\begin{bmatrix}s\\k\end{bmatrix}_{q_i}[\frac{1}{vq_i^{-2}-1}-\frac{q_i^{2k}}{vq_i^{2k-2}-1}]$

\noindent $\mu_k=(-1)^{k-1}w_1^{s-1}\begin{bmatrix}s\\k-1\end{bmatrix}_{q_i}\underset{k'=k...s}{\sum}(v-q_i^2)(vq_i^2-q_i^2)...(vq_i^{2(k'-2)}-q_i^2)(vq_i^{2k'}-1)...(vq_i^{2(s-2)+2}-1)
\\=Q(-1)^{k-1}\begin{bmatrix}s\\k-1\end{bmatrix}_{q_i}[\frac{q_i^{2k-2}}{vq_i^{2k-4}-1}-\frac{q_i^{2s}}{vq_i^{2s-2}-1}]$

\noindent As $\alpha_k=q^{1-s}\lambda_k+\mu_k$ and $\beta_k=-(\lambda_k+q_i^{1-s})/(q_i^{k-2}w_2)$, we have:
$$\alpha_k=Q\frac{(-1)^k\begin{bmatrix}s\\k\end{bmatrix}_{q_i}[k]_{q_i}(q_i-q_i^{-1})((q_i^{k+1-s}-q_i^{s+k-1})+v(q_i^{s+k-3}+q_i^{s+3k-3}-q_i^{3k-3-s}-q_i^{s+k-1}))}{(vq_i^{-2}-1)(vq_i^{2k-2}-1)(vq_i^{2k-4}-1)(vq_i^{2s-2}-1)}$$
$$\beta_k=Q\frac{(-1)^k\begin{bmatrix}s\\k\end{bmatrix}_{q_i}[k]_{q_i}(q_i-q_i^{-1})((q_i^k+q_i^{2s-k+2}-q_i^{-k+2}-q_i^{k+2})+v(-q_i^{k+2s-2}+q_i^k))}{(vq_i^{-2}-1)(vq_i^{2k-2}-1)(vq_i^{2k-4}-1)(vq_i^{2s-2}-1)}$$
In particular $(q_i^2-q_i^{2-2s})f_{+,k-1}^{(k)}(Q(-1)^k\begin{bmatrix}s\\k\end{bmatrix}_{q_i}[k]_{q_i}(q_i-q_i^{-1}))^{-1}$ is:
$$\frac{(q_i^{k+1-s}+q_i^{s-k+3}-q_i^{s+k+1}-q_i^{3-k-s})+v(q_i^{s+3k-1}+q_i^{k+1-s}-q_i^{3k-1-s}-q_i^{s+k+1})}{(vq_i^{-2}-1)(vq_i^{2k-2}-1)(vq_i^{2k-4}-1)(vq_i^{2s-2}-1)}$$
and we get formula \ref{trouvun} for $k$. Moreover $(q_i^2-q_i^{2-2s})f_{+,s-1}^{(k)}(Q(-1)^{k}\begin{bmatrix}s\\k\end{bmatrix}_{q_i}[k]_{q_i}(q_i-q_i^{-1}))^{-1}$ is:
$$\frac{(q_i^{k+2-2s}+q_i^{2s-k+2}-q_i^{-k+2}-q_i^{k+2})+v(q_i^{k-2}+q_i^{3k-2}-q_i^{3k-2s-2}-q_i^{k+2s-2})}{(vq_i^{-2}-1)(vq_i^{2k-2}-1)(vq_i^{2k-4}-1)(vq_i^{2s-2}-1)}$$
and we get formula \ref{trouvdeux} for $k+1$. 

\noindent So it remains to prove formula \ref{trouvdeux} with $k=2$ and formula \ref{trouvun} with $k=s$. 
$$P_+^{(1)}=(w_2-q_i^{(1-s)}z)f_{+,s-1}^{(1)}=-[s]_{q_i}(zq_i^{1-s}-w_1)(q_i^2w_2-w_1)...(q_i^{2s-2}w_2-w_1)$$
$$+(z-w_1q_i^{1-s})\underset{k'=1...s}{\sum}q_i^{2k'-2}(q_i^{-2}w_2-w_1)...(q_i^{2k'-6}w_2-w_1)(q_i^{2k'}w_2-w_1)...(q_i^{2s-2}w_2-w_1)$$
$$\Rightarrow f_{+,s-1}^{(1)}=-q_i^{1-s}Q[\frac{-[s]_{q_i}q_i^{1-s}(q_i^2-1)}{(vq_i^{-2}-1)(v-1)}+\underset{k'=1...s}{\sum}\frac{q_i^{2k'-2}}{(q_i^{2k'-4}v-1)(q_i^{2k'-2}v-1)}]$$
And so we have for $f_{+,s-1}^{(1)}(q_i^2-q_i^{2-2s})(-Q[s]_{q_i}(q_i-q_i^{-1}))^{-1}$:
$$\frac{q_i+q_i^3-q_i^{2s+1}-q_i{-2s+3}}{(vq_i^{-2}-1)(v-1)(vq_i^{2s-2}v-1)}$$
that it to say the formula \ref{trouvdeux} with $k=2$.
$$P_+^{(s)}=(z-q_i^{(1-s)}q_i^{2(s-2)}w_2)f_{+,s-1}^{(s)}=(-1)^{s-1}[s]_{q_i}(z-w_1q_i^{1-s})q_i^{2(s-1)}(q_i^{-2}w_2-w_1)...(q_i^{2s-6}w_2-w_1)$$
$$+(-1)^s(zq_i^{1-s}-w_1)\underset{k'=1...s}{\sum}q_i^{2k'-2}(q_i^{-2}w_2-w_1)...(q_i^{2k'-6}w_2-w_1)(q_i^{2k'}w_2-w_1)...(q_i^{2s-2}w_2-w_1)$$
$$\Rightarrow f_{+,s-1}^{(s)}=Q[\frac{(-1)^{s-1}[s]_{q_i}q_i^{2(s-1)}(q_i^2-1)}{(vq_i^{2s-4}-1)(vq_i^{2s-2}-1)}+(-1)^sq_i^{1-s}\underset{k'=1...s}{\sum}\frac{q_i^{2k'-2}}{(q_i^{2k'-4}v-1)(q_i^{2k'-2}v-1)}]$$
And so we have for $f_{+,s-1}^{(s)}(q_i^2-q_i^{2-2s})((-1)^sQ [s]_{q_i}(q_i-q_i^{-1}))^{-1}$:
$$\frac{q_i^{2s+1}+q_i^{3-2s}-q_i-q_i^3}{(vq_i^{-2}-1)(vq_i^{2s-4}-1)(vq_i^{2s-2}v-1)}$$
that it to say the formula \ref{trouvun} with $k=s$.\qed

\subsubsection{Proof of theorem \ref{dectrian}}\label{finpreuve} The algebras $\U_q^l(\hat{\Glie}), \tilde{\U}_q(\hat{\Glie}), \U_q^{l, \pm}(\hat{\Glie})$ are defined in section \ref{proofdecomp}. Let $\tilde{\U}_q^{\pm}(\hat{\Glie})\subset\tilde{\U}_q(\hat{\Glie})$ be the subalgebra generated by the $x_{i,r}^{\pm}$. Let $\tau_{\pm}$ be the two-sided ideal of $\U_q^{l,\pm}(\hat{\Glie})$ generated by the left terms of relations (\ref{equaun}) (with the $x_{i,r}^{\pm}$). 

\begin{lem}\label{stepun} We have $\tau_+\U_q^l(\hat{\Glie})\subset \U_q^l(\hat{\Glie})\tau_+$ and $\U_q^l(\hat{\Glie})\tau_-\subset \tau_-\U_q^l(\hat{\Glie})$. In particular $(\tilde{\U}_q^-(\hat{\Glie}),\U_q(\hat{\Hlie}),\tilde{\U}_q^+(\hat{\Glie}))$ is a triangular decomposition of $\tilde{\U}_q(\hat{\Glie})$.\end{lem}

\demo First $\tau_+\U_q(\hat{\Hlie})\subset \U_q(\hat{\Hlie})\tau_+$, $\U_q(\hat{\Hlie})\tau_-\subset\tau_-\U_q(\hat{\Hlie})$ are direct consequences of relations (\ref{actcartc}), (\ref{actcartplus}), (\ref{actcartmoins}). We have also (we use relations (\ref{partact}) and (\ref{actcartplus}), (\ref{actcartmoins})):
$$[(w-q^{\pm B_{i,j}}z)x_i^{\pm}(z)x_j^{\pm}(w)-(q^{\pm B_{i,j}}w-z)x_j^{\pm}(w)x_i^{\pm}(z),x_k^{\mp}(u)]$$
$$=(w-q^{\pm B_{i,j}}z)x_i^{\pm}(z)[x_j^{\pm}(w),x_k^{\mp}(u)]
-(q^{\pm B_{i,j}}w-z)[x_j^{\pm}(w),x_k^{\mp}(u)]x_i^{\pm}(z)$$
$$-(q^{\pm B_{i,j}}w-z)x_j^{\pm}(w)[x_i^{\pm}(z),x_k^{\mp}(u)]
+(w-q^{\pm B_{i,j}}z)[x_i^{\pm}(z),x_k^{\mp}(u)]x_j^{\pm}(w)=0$$
and so $\tau_+\U_q^{l,-}(\hat{\Glie})\subset \U_q^l(\hat{\Glie})\tau_+$, $\U_q^{l,+}(\hat{\Glie})\tau_-\subset \tau_-\U_q^l(\hat{\Glie})$.

\noindent The last point follows from $\tilde{\U}_q(\hat{\Glie})=\U_q^l(\hat{\Glie})/(\U_q^l(\hat{\Glie}).(\tau_++\tau_-).\U_q^l(\hat{\Glie}))$, the triangular decomposition of $\U_q^l(\hat{\Glie})$ and lemma \ref{gendec}.\qed

\begin{lem}\label{stepdeux} Let $i\neq j$, $s=1-C_{i,j}$ $\mu=1\text{ or }\mu=-1$. We have in $\tilde{\U}_q(\hat{\Glie})$ :
\begin{equation}\label{tildeun}\underset{\pi\in \Sigma_s}{\sum}\underset{k=0..s}{\sum}(-1)^k\begin{bmatrix}s\\k\end{bmatrix}_{q_i}x_i^{\pm}(w_{\pi(1)})...x_i^{\pm}(w_{\pi(k)})\phi^{\mu}_j(z)x_i^{\pm}(w_{\pi(k+1)})...x_i^{\pm}(w_{\pi(s)})=0\end{equation}
\begin{equation}\label{tildedeux}\underset{\pi\in \Sigma_s}{\sum}\underset{k=0..s}{\sum}(-1)^k\begin{bmatrix}s\\k\end{bmatrix}_{q_i}\xi_i(w_{\pi(1)})...\xi_i(w_{\pi(k)})x_j^{\pm}(z)\xi_i(w_{\pi(k+1)})...\xi_i(w_{\pi(s)})=0\end{equation}
where $\xi_i(w_p)=x_i^{\pm}(w_p)$ if $p\neq 1$ and $\xi_i(w_1)=\phi_i^{\mu}(w_1)$.\end{lem}

\noindent Remark : in particular if we multiply the equation (\ref{tildeun}) by $(\underset{r=1...s}{\prod}(w_r-q_i^{s-1}z))(\underset{1\leq r'<r\leq s}{\prod}(w_r-q_i^2w_{r'}))$ and we project it on $x_i^+(w_1)...x_i^+(w_s)\phi_j^+(z)$ (we can use the relations (\ref{plusmoinsc}) thanks to the multiplied polynomial), we get the combinatorial identity discovered by Jing in \cite{jin}, which was also proved in a combinatorial way in \cite{dj} : for $\pi\in\Sigma_s$ denote by $\epsilon(\pi)\in\{1,-1\}$ the signature of $\pi$ (we have replaced $z\mapsto z^{-1}$, $w_{k'}\mapsto w_{k'}^{-1}$ to get the formula in the same form as in \cite{jin}):
$$0=\underset{\pi\in\Sigma_s}{\sum}\epsilon(\pi)\underset{k=0...s}{\sum}\begin{bmatrix}s\\k\end{bmatrix}_q(z-q^{s-1}w_{\pi(1)})...(z-q^{s-1}w_{\pi(k)})$$
$$(w_{\pi(k+1)}-q^{s-1}z)...(w_{\pi(s)}-q^{s-1}z)\underset{1\leq r<r'\leq s}{\prod}(w_{\pi(r)}-q^2w_{\pi(r')})$$

\demo First we prove the equation (\ref{tildeun}) with $\mu=1$ ($\mu=-1$ is analog). The left term is (relations (\ref{actcartplus})):
$$\frac{\phi^+_j(z)}{(w_1 q^{\pm B_{i,j}}-z)...(w_s q^{\pm B_{i,j}}-z)}\underset{\pi\in \Sigma_s}{\sum}P_{\pm}(w_{\pi(1)},...,w_{\pi(s)},z)x_i^{\pm}(w_{\pi(1)})...x_i^{\pm}(w_{\pi(s)})$$
that is to say (see lemma \ref{techun}): 
$$\underset{\pi\in \Sigma_s}{\sum}\underset{1\leq r\leq s-1}{\sum}(w_{\pi(r+1)}-q_i^{\pm 2}w_{\pi(r)})f_{r ,\pm}(w_{\pi(1)},...,w_{\pi(r-1)},w_{\pi(r+2)},...,w_{\pi(s)},z)x_i^{\pm}(w_{\pi(1)})...x_i^{\pm}(w_{\pi(s)})$$
For each $r$, we put together the $\pi,\pi'\in\Sigma_s$ such that $\pi(r)=\pi'(r+1)$, $\pi(r+1)=\pi'(r)$, and $\pi(r'')=\pi'(r'')$ for all $r''\neq r, r+1$. So we get a sum of terms:
$$f_{r ,\pm}(w_{\pi(1)},...,w_{\pi(r-1)},w_{\pi(r+2)},...,w_{\pi(s)},z)x_i^{\pm}(w_{\pi(1)})...x_i^{\pm}(w_{\pi(r-1)})A^{\pm}_{\{\pi(r),\pi(r+1)\}}x_i^{\pm}(w_{\pi(r+2})...x_i^{\pm}(w_{\pi(s)})$$
$$\text{where } A^{\pm}_{\{k,k'\}}=(w_{k}-q_i^{\pm 2} w_{k'})x_i^{\pm}(w_{k'})x_i^{\pm}(w_{k})+(w_{k'}-q_i^{\pm 2} w_{k})x_i^{\pm}(w_{k})x_i^{\pm}(w_{k'})$$ 
But $A^{\pm}_{\{k,k'\}}=0$ in $\tilde{\U}_q(\Glie)$.

\noindent Let us prove the equation (\ref{tildedeux}) with $\mu=1$ ($\mu=-1$ is analog). The left term is :
$$\frac{\phi^+_i(w_1)}{(w_2 q_i^{\pm 2} -w_1)...(w_s q_i^{\pm 2}-w_1)(zq_i^{\pm (1-s)}-w_1)}$$
$$\underset{\pi\in \Sigma_{s-1}, k=1 ... s}{\sum}P_{\pm}^{(k)}(w_1,w_{\pi(2)},...,w_{\pi(s)},z)x_i^{\pm}(w_{\pi(2)})...x_i^{\pm}(w_{\pi(k)})x_j^{\pm}(z)x_i^{\pm}(w_{\pi(k+1)})...x_i^{\pm}(w_{\pi(s)})$$
where $\Sigma_{s-1}$ acts on $\{2,...,s\}$. With the help of lemma \ref{techdeux} and in analogy to the previous case, for each $1\leq k\leq s$ each $r\neq k$, we put together the $\pi,\pi'\in\Sigma_s$ such that $\pi(r)=\pi'(r+1)$, $\pi(r+1)=\pi'(r)$, and $\pi(r'')=\pi'(r'')$ for all $r''\neq r, r+1$. So the terms with polynomials $f_{\pm,k'}^{(k)}$ with $k'\neq s, k-1$ are erased. We get : $\frac{\phi^+_i(w_1)}{(w_2 q_i^{\pm 2} -w_1)...(w_s q_i^{\pm 2}-w_1)(zq_i^{\pm (1-s)}-w_1)}
\\ \underset{\pi\in \Sigma_{s-1}, k=1 ... s}{\sum} ((z-q_i^{\pm (1-s)}w_{\pi(k)})f_{\pm ,k-1}^{(k)}+(w_{\pi(k+1)}-q_i^{\pm (1-s)}z)f_{\pm, s-1}^{(k)})
\\x_i^{\pm}(w_{\pi(2)})...x_i^{\pm}(w_{\pi(k)})x_j^{\pm}(z)x_i^{\pm}(w_{\pi(k+1)})...x_i^{\pm}(w_{\pi(s)})$

\noindent But this last sum is equal to :
$$\underset{\pi\in \Sigma_{s-1}, k=2 ... s}{\sum} (z-q_i^{\pm (1-s)}w_{\pi(k)})(f_{\pm ,k-1}^{(k)}-f_{\pm, s-1}^{(k-1)})x_i^{\pm}(w_{\pi(2)})...x_i^{\pm}(w_{\pi(k)})x_j^{\pm}(z)x_i^{\pm}(w_{\pi(k+1)})...x_i^{\pm}(w_{\pi(s)})$$
where we can replace $(z-q_i^{\pm (1-s)}w_{\pi(k)})x_i^{\pm}(w_{\pi(k)})x_j^{\pm}(z)$ by $(-w_{\pi(k)}+q_i^{\pm (1-s)}z)x_j^{\pm}(z)x_i^{\pm}(w_{\pi(k)})$ (relations (\ref{plusmoinsc}) in $\tilde{\U}_q(\hat{\Glie})$). As in the previous cases it follows from lemma \ref{techtrois} that this term is equal to $0$.\qed

\noindent Let $\tilde{\tau}_{\pm}$ be the two-sided ideal of $\tilde{\U}_q^{\pm}(\hat{\Glie})$ generated by the left terms of relations (\ref{equadeux}) with the $x_{i,r}^{\pm}$. 

\begin{lem}\label{steptrois} We have $\tilde{\tau}_+\tilde{\U}_q(\hat{\Glie})\subset \tilde{\U}_q(\hat{\Glie})\tilde{\tau}_+$ and $\tilde{\U}_q(\hat{\Glie})\tilde{\tau}_-\subset \tilde{\tau}_-\tilde{\U}_q(\hat{\Glie})$. \end{lem}

\noindent In particular as $\U_q(\hat{\Glie})=\tilde{\U}_q^l(\hat{\Glie})/(\tilde{\U}_q^l(\hat{\Glie}).(\tilde{\tau}_++\tilde{\tau}_-).\tilde{\U}_q^l(\hat{\Glie}))$ the result of theorem \ref{dectrian} follows from lemma \ref{gendec} and the triangular decomposition of $\tilde{\U}_q(\hat{\Glie})$ proved in lemma \ref{stepun}.

\demo First $\tilde{\tau}_+\U_q(\hat{\Hlie})\subset\U_q(\hat{\Hlie})\tau_+$, $\U_q(\hat{\Hlie})\tilde{\tau}_-\subset\tilde{\tau}_-\U_q(\hat{\Hlie})$ are direct consequences of relations (\ref{actcartc}), (\ref{actcartplus}), (\ref{actcartmoins}). Let us show that :
\begin{equation}\label{aijl}[\underset{\pi\in \Sigma_s}{\sum}\underset{k=0..s}{\sum}(-1)^k\begin{bmatrix}s\\k\end{bmatrix}_{q_i}x_i^{\pm}(w_{\pi(1)})...x_i^{\pm}(w_{\pi(k)})x_j^{\pm}(z)x_i^{\pm}(w_{\pi(k+1)})...x_i^{\pm}(w_{\pi(s)}),x_l^{\mp}(u)]=0\end{equation}
where $i,j,l\in I$, $i\neq j$. If $l\neq j$ and $l\neq i$ the equation (\ref{aijl}) follows from relations (\ref{partact})). If $l=j$, the equation (\ref{aijl}) follows from the identity (\ref{tildeun}) of lemma \ref{stepdeux} because the left term is :
$$\underset{\pi\in \Sigma_s}{\sum}\underset{k=0..s}{\sum}(-1)^k\begin{bmatrix}s\\k\end{bmatrix}_{q_i}x_i^{\pm}(w_{\pi(1)})...x_i^{\pm}(w_{\pi(k)})(\delta(\frac{z}{u})\phi^{\pm}_j(z)-\delta(\frac{z}{u})\phi^{\mp}_j(z))x_i^{\pm}(w_{\pi(k+1)})...x_i^{\pm}(w_{\pi(s)})$$

\noindent If $l=i$, the equation (\ref{aijl}) follows from the identity (\ref{tildedeux}) of lemma \ref{stepdeux} because the left term is :
$$\underset{\pi\in \Sigma_s}{\sum}\underset{k=0..s}{\sum}(-1)^k\begin{bmatrix}s\\k\end{bmatrix}_{q_i}(\underset{k'=1...k}{\sum}x_i^{\pm}(w_{\pi(1)})...x_i^{\pm}(w_{\pi(k'-1)})\delta(\frac{w_{k'}}{u})(\phi^{\pm}_i(w_{k'})-\phi^{\mp}_i(w_{k'}))$$
$$x_i^{\pm}(w_{\pi(k'+1)})...x_i^{\pm}(w_{\pi(k)})x_j^{\pm}(z)x_i^{\pm}(w_{\pi(k+1)})...x_i^{\pm}(w_{\pi(s)})$$
$$+\underset{k'=k+1 ... s}{\sum}x_i^{\pm}(w_{\pi(1)})...x_i^{\pm}(w_{\pi(k)})x_j^{\pm}(z)x_i^{\pm}(w_{\pi(k+1)})...x_i^{\pm}(w_{\pi(k'-1)})$$
$$\delta(\frac{w_{k'}}{u})(\phi^{\pm}_i(w_{k'})-\phi^{\mp}_i(w_{k'}))x_i^{\pm}(w_{\pi(k'+1)})...x_i^{\pm}(w_{\pi(s)}))$$
So we have proved the equation (\ref{aijl}) and in particular $\tilde{\tau}_+\tilde{\U}_q^-(\hat{\Glie})\subset \tilde{\U}_q(\hat{\Glie})\tilde{\tau}_+$, $\tilde{\U}_q^+(\hat{\Glie})\tilde{\tau}_-\subset \tilde{\tau}_-\tilde{\U}_q(\hat{\Glie})$.\qed

\section{Integrable representations and category $\mathcal{O}(\U_q(\hat{\Glie}))$}\label{int}

In this section we study highest weight representations of $\U_q(\hat{\Glie})$. In particular theorem \ref{simpint} is a generalization of a result of Chari-Pressley about integrable representations.

\subsection{Reminder: integrable representations of quantum Kac-Moody algebras}\label{rem} In this section we review some known properties of integrable representations of $\U_q(\Glie)$.

\noindent For $V$ a $\U_q(\Hlie)$-module and $\omega\in \Hlie^*$ we denote by $V_{\omega}$ the weight space of weight $\omega$:
$$V_{\omega}=\{v\in V/\forall h\in \Hlie, k_h.v=q^{\omega(h)}v\}$$
In particular for $v\in V_{\omega}$ we have $k_i.v=q_i^{\omega(\alpha_i^{\vee})}v$ and for $i\in I$ we have $x_i^{\pm}.V_{\omega}\subset V_{\omega \pm \alpha_i}$.

\noindent We say that $V$ is $\U_q(\Hlie)$-diagonalizable if $V=\underset{\omega\in \Hlie^*}{\bigoplus}V_{\omega}$ (in particular $V$ is of type $1$).

\begin{defi} A $\U_q(\Glie)$-module $V$ is said to be integrable if $V$ is $\U_q(\Hlie)$-diagonalizable, $\forall\omega\in \Hlie^*$, $V_{\omega}$ is finite dimensional, and for $\mu\in \Hlie^*$, $i\in I$ there is $R\geq 0$ such that $r\geq R \Rightarrow V_{\mu\pm r\alpha_i}=\{0\}$.\end{defi}

\noindent In particular for all $v\in V$ there is $m_v\geq 0$ such that for all $i\in I$, $m\geq m_v$, $(x_i^+)^m.v=(x_i^-)^m.v=0$, and $U_i.v$ is finite dimensional.

\begin{defi} A $\U_q(\Glie)$-module $V$ is said to be of highest weight $\omega\in \Hlie^*$ if there is $v\in V_{\omega}$ such that $V$ is generated by $v$ and $\forall i\in I$, $x_i^+.v=0$.
\end{defi}

\noindent In particular $V=\U_q(\Glie)^-.v$ (theorem \ref{triansc}), $V$ is $\U_q(\Hlie)$-diagonalizable, and $V=\underset{\lambda \leq \omega}{\bigoplus}V_{\lambda}$. We have (see \cite{lu}):

\begin{thm} For any $\omega\in \Hlie^*$ there is a unique up to isomorphism simple highest weight module $L(\omega)$ of highest weight $\omega$. The highest weight module $L(\omega)$ is integrable if and only $\omega\in P^+$.\end{thm}

\subsection{Integrable representations of quantum affinizations} In this section we generalize results of Chari-Pressley \cite{Cha, Cha2} to all quantum affinizations. 

\subsubsection{$l$-highest weight modules}\label{needtrian}

We give the following notion of $l$-weight :

\begin{defi} A couple $(\lambda,\Psi)$ such that $\lambda\in \Hlie^*$, $\Psi=(\Psi_{i,\pm m}^{\pm})_{i\in I, m\geq 0}$, $\Psi_{i,\pm m}^{\pm}\in\CC$, $\Psi_{i,0}^{\pm}=q_i^{\pm \lambda(\alpha_i^{\vee})}$ is called a $l$-weight. \end{defi}
\noindent The condition $\Psi_{i,0}^{\pm}=q_i^{\pm \lambda(\alpha_i^{\vee})}$ is a compatibility condition which comes from $\phi_{i,0}^{\pm}=k_i^{\pm}$.

\noindent We denote by $P_l$ the set of $l$-weights. Note that in the finite case $\lambda$ is uniquely determined by $\Psi$ because $\lambda=\underset{i\in I}{\sum}\lambda(\alpha_i^{\vee})\Lambda_i$. Analogs of those $l$-weights were also used in \cite{mi} for toroidal $\hat{sl_n}$-cases.

\begin{defi} A $\U_q(\hat{\Glie})$-module $V$ is said to be of $l$-highest weight $(\lambda, \Psi)\in P_l$ if there is $v\in V$ such that ($i\in I, r\in\ZZ, m\geq 0, h\in \Hlie$):
$$x_{i,r}^+.v=0\text{ , }V=\U_q(\hat{\Glie}).v\text{ , }\phi_{i,\pm m}^{\pm}.v=\Psi_{i,\pm m}^{\pm}v\text{ , }k_h.v=q^{\lambda(h)}.v$$
\end{defi}

\noindent In particular $\U_q(\hat{\Glie})^-.v=V$ (theorem \ref{dectrian}), $V$ is $\U_q(\Hlie)$-diagonalizable and $V=\underset{\lambda\leq \omega}{\bigoplus}V_{\lambda}$. Note that the $l$-weight $(\lambda, \Psi)\in P_l$ is uniquely determined by $V$. It is called the $l$-highest weight of $V$.

\noindent The notion of $l$-highest weight is different from the notion of highest weight for quantum affine algebras. The term ``pseudo highest weight'' is also used in the literature. 

\noindent Example : for any $(\lambda, \Psi)\in P_l$, define the Verma module $M(\lambda, \Psi)$ as the quotient of $\U_q(\hat{\Glie})$ by the left ideal generated by the $x_{i,r}^+$ ($i\in I, r\in\ZZ$), $k_h-q^{\lambda(h)}$ ($h\in \Hlie$), $\phi_{i,\pm m}^{\pm}-\Psi_{i,\pm m}^{\pm}$ ($i\in I, m\geq 0$). It follows from theorem \ref{dectrian} that $M(\lambda, \Psi)$ is a free $\U_q^-(\hat{\Glie})$-module of rank $1$. In particular it is non trivial and it is a $l$-highest weight module of highest weight $(\lambda, \Psi)$. Moreover it has a unique proper submodule (mimic the classical argument in \cite{kac}), and :

\begin{prop}\label{simple} For any $(\lambda, \Psi)\in P_l$ there is a unique up to isomorphism simple $l$-highest weight module $L(\lambda, \Psi)$ of $l$-highest weight $(\lambda, \Psi)$.\end{prop}

\subsubsection{Integrable $\U_q(\hat{\Glie})$-modules}

\begin{defi}\label{defint} A $\U_q(\hat{\Glie})$-module $V$ is said to be integrable if $V$ is integrable as a $\U_q(\Glie)$-module.\end{defi}

\noindent Note that in the case of a quantum affine algebra, the two notions of integrability do not coincide. Throughout the paper only the notion of integrability of definition \ref{defint} is used.

\noindent For $i\in I, r\in \ZZ$ and $\omega\in \Hlie^*$ we have $x_{i,r}^{\pm}.V_{\omega}\subset V_{\omega \pm \alpha_i}$. So if $V$ is integrable, for all $v\in V$, $\hat{U}_i.v$ is finite dimensional and there is $m_0\geq 1$ such that for all $i\in I$, $r\in\ZZ$, $m\geq m_0\Rightarrow (x_{i,r}^+)^m.v=(x_{i,r}^-)^m.v=0$.

\begin{defi} The set $P_l^+$ of dominant $l$-weights is the set of $(\lambda,\Psi)\in P_l$ such that there exist (Drinfel'd)-polynomials $P_i(z)\in\CC[z]$ ($i\in I$) of constant term $1$ such that in $\CC[[z]]$ (resp. in $\CC[[z^{-1}]]$):
$$\underset{m\geq 0}{\sum} \Psi_{i,\pm m}^{\pm} z^{\pm m}=q_i^{\text{deg}(P_i)}\frac{P_i(zq_i^{-1})}{P_i(zq_i)}$$
\end{defi}

\noindent In particular for all $i\in I$, $\lambda(\alpha_i^{\vee})=\text{deg}(P_i)\geq 0$ and so $\lambda\in P^+$ is a dominant weight.

\begin{thm}\label{simpint} For $(\lambda, \Psi)\in P_l$, $L(\lambda, \Psi)$ is integrable if and only $(\lambda, \Psi)\in P_l^+$.\end{thm}

\noindent If $\Glie$ is finite (case of a quantum affine algebra) it is a result of Chari-Pressley in \cite{Cha} (if part) and in \cite{Cha2} (only if part). Moreover in this case the integrable $L(\lambda,\Psi)$ are finite dimensional. If $\Glie$ is symmetric the result is geometrically proved by Nakajima in \cite{Naams}. If $C$ is of type $A_n^{(1)}$ (toroidal $\hat{sl_n}$-case) the result is algebraically proved by Miki in \cite{mi}.

\noindent For the general case we propose a proof similar to the proof given by Chari-Pressley in the finite case. For $\lambda\in \Hlie^*$ denote $D(\lambda )=\{\omega\in \Hlie^*/\omega\leq\lambda\}$.

\demo The proof uses the result for $\U_q(\hat{sl_2})$ which is proved in \cite{Cha0, Cha}.

\noindent First suppose that $L=L(\lambda, \Psi)$ is integrable and for $i\in I$ let $L_i$ be the $\hat{U}_i$-submodule of $L$ generated by the highest weight vector $v$. It is a $l$-highest weight $\U_{q_i}(\hat{sl_2})$-module of highest weight $(\lambda(\alpha_i^{\vee}), \Psi_i^{\pm})$. As $L$ is integrable, $L_i$ is finite dimensional. So the result for $\U_{q_i}(\hat{sl_2})$ gives $P_i(z)\in\CC[z]$ such that:
$$\underset{m\geq 0}{\sum}\Psi_{i,\pm m}^{\pm}z^{\pm m}=q_i^{\text{deg}(P_i)}\frac{P_i(zq_i^{-1})}{P_i(zq_i)}\text{ , }\lambda(\alpha_i^{\vee})=\text{deg}(P_i)\geq 0$$

\noindent Now we prove that $L=L(\lambda, \Psi)=\U_q(\hat{\Glie}).v$ is integrable where $(\lambda,\Psi)\in P_l^+$. It suffices to prove that:

(1) For all $\mu\leq \lambda$, if $L_{\mu}\neq \{0\}$ then there exists $M>0$ such that $m>M\Rightarrow L_{\mu - m\alpha_i}=L_{\mu + m\alpha_i}=0$ for all $i\in I$.

(2) For all $\mu\leq \lambda$, $\text{dim}(L_{\mu})<\infty$.

\noindent The proof goes roughly as in section 5 of \cite{Cha}, with the following modifications : 

For (1) : the existence of $M$ for $L_{\mu + m\alpha_i}=0$ is clear because the weights of $L$ are in $D(\lambda)$. Put $r^{\vee}=\text{max}\{-C_{i,j}/i\neq j\}$. In particular if $C$ is finite, we have $r^{\vee}\leq 3$. First we prove that for $m>0$, the space $L_{\mu-m\alpha_i}$ is spanned by vectors of the form $X_1^-x_{i_1,k_1}^-...X_h^-x_{i_h,k_h}^-X_{h+1}^-.v$ where $\lambda-\mu=\alpha_{i_1}+...+\alpha_{i_h}$, $k_1,...,k_h\in\ZZ$, $X_p^-$ is of the form $X_p^-=x_{i,l_{1,p}}^-...x_{i,l_{m_p,p}}^-$ where $m_1+...+m_{h+1}=m$ and $m_1,...,m_h\leq r^{\vee}$ (which is the crucial condition). It is proved by induction on $h$ (see \cite{Cha} section 5, (e)) with the help of the relations (\ref{equadeux}). Note that in \cite{Cha} $r^{\vee}=3$. Now it suffices to prove that $\hat{U}_i.v$ is finite dimensional : indeed if $m>r^{\vee}h+\text{dim}(\hat{U}_i.v)$ we have $m_{h+1}>\text{dim}(\hat{U}_i.v)$ and $X_{h+1}^-.v=0$. It is shown exactly as in lemma 2.3 of \cite{Chah} that $\hat{U}_i.v$ is irreducible as $\hat{U}_i$-module, and so is finite dimensional.

For (2) : let us write $\lambda-\mu=\alpha_{i_1}+...+\alpha_{i_h}$. The result is proved by induction on $h$. We have seen that $\hat{U}_i.v$ is finite dimensional. The induction is shown exactly as in \cite{Cha} (section 5. (b)) by considering the $L_{\lambda-\mu+\alpha_{i_j}}$ and with the help of relation (\ref{equaun}).\qed

\subsection{Category $\mathcal{O}(\U_q(\hat{\Glie}))$} In the following by subcategory we mean full subcategory. 

\begin{defi} A $\U_q(\Hlie)$-module $V$ is said to be in the category $\mathcal{O}(\U_q(\Hlie))$ if:

i) $V$ is $\U_q(\Hlie)$-diagonalizable

ii) for all $\omega\in \Hlie^*$, $\text{dim}(V_{\omega})<\infty$

iii) there is a finite number of element $\lambda_1,...,\lambda_s\in \Hlie^*$ such that the weights of $V$ are in $\underset{j=1...s}{\bigcup}D(\lambda_j)$

\noindent A $\U_q(\Glie)$-module (resp. a $\U_q(\hat{\Glie})$-module) is said to be in the category $\mathcal{O}(\U_q(\Glie))$ (resp. $\mathcal{O}(\U_q(\hat{\Glie}))$) if it is in the category $\mathcal{O}(\U_q(\Hlie))$ as a $\U_q(\Hlie)$-module.\end{defi}

\noindent  In particular we have a restriction functor $\text{res}:\mathcal{O}(\U_q(\hat{\Glie}))\rightarrow\mathcal{O}(\U_q(\Glie))$.

\noindent For example a highest weight $\U_q(\Glie)$-module is in the category $\mathcal{O}(\U_q(\Glie))$ and the product $\otimes$ is well-defined on $\mathcal{O}(\U_q(\Glie))$. An integrable $l$-highest weight module is in the category $\mathcal{O}(\U_q(\hat{\Glie}))$. But in general a $l$-highest weight module is not in the category $\mathcal{O}(\U_q(\hat{\Glie}))$, indeed ($\CC_r[z]$ is the space of polynomials of degree lower that $r$):

\begin{lem} Consider a $l$-weight $(\omega,\Psi)\in P_l$ and $i\in I$. If $\text{dim}(L(\omega,\Psi)_{\omega-\alpha_i})=r\in \NN$ then there is $P(z)\in\CC_r[z]$ such that $P(z)\Psi_i(z)=0$ where $\Psi_i(z)=\underset{r\geq 0}{\sum}(\Psi_{i,r}^+z^r-\Psi_{i,-r}^-z^{-r})$.\end{lem}

\noindent In particular the existence of a $P(z)\in\CC[z]$ such that $P(z)\Psi_i(z)=0$ for all $i\in I$ is a necessary condition for $L(\omega,\Psi)\in\mathcal{O}(\U_q(\hat{\Glie}))$.

\demo Let $v_0,v_1,...,v_r\in L(\omega,\Psi)$ such that:
$$L(\omega,\Psi)_{\omega}=\CC v_0\text{ , }L(\omega,\Psi)_{\omega-\alpha_i}=\CC v_1\oplus ...\oplus\CC v_r$$
For $m\in\ZZ$ let $\Psi_{i,m}=\Psi_{i,m}^+ -\Psi_{i,m}^{-}$. As $x_{i.m}^+.v_0=0$, we have:
$$x_{i,m}^+x_{i,m'}^-.v_0=\frac{1}{q_i-q_i^{-1}}\Psi_{i,m+m'}v_0$$
As $x_{i,m}^-.v_0\in L(\omega,\Psi)_{\omega -\alpha_i}$ and $x_{i,m}^+.v_j\in L(\omega,\Psi)_{\omega}$, there are $\lambda_m^j,\mu_m^j\in\CC$ ($m\in\ZZ,1\leq j\leq r$) such that:
$$x_{i,m}^-.v_0=\lambda_m^1v_1+...+\lambda_m^rv_r\text{ , }x_{i,m}^+.v_j=\mu_m^jv_0$$
In particular we have: $\Psi_{i,m+m'}=(q_i-q_i^{-1})\underset{j=1...r}{\sum}\lambda_{m'}^j\mu_m^j$. We set $\lambda^j(z)=\underset{m'\in\ZZ}{\sum}\lambda^j_{m'}z^{m'}$, $\Psi_i(z)=\underset{r\geq 0}{\sum}\Psi_{i,r}^+z^r-\Psi_{i,-r}^-z^{-r}$ and we have :
$$z^{-m}\Psi_i(z)=(q_i-q_i^{-1})\underset{j=1...r}{\sum}\mu_m^j\lambda^j(z)$$
So the $\{\Psi_i(z),z\Psi_i(z),...,z^r\Psi_i(z)\}$ are not linearly independent.\qed

\section{$q$-characters}\label{qcar} For a quantum Kac-Moody algebra, one can define a character morphism as in the classical case. For quantum affine algebras a more precise morphism, called morphism of $q$-characters, was introduced by Frenkel-Reshetikhin \cite{Fre} (in particular to distinguish finite dimensional representations). In this section we generalize the construction of $q$-characters to quantum affinizations. The technical point is to add terms $k_{\lambda}$ ($\lambda\in\Hlie^*$) to make it well-defined in the general case. We prove a symmetry property of $q$-characters that generalizes a result of Frenkel-Mukhin : the image of $\chi_q$ is the intersection of the kernels of screening operators (theorem \ref{sym}).

\subsection{Reminder: classical character} Let $\U_q(\Glie)$ be a quantum Kac-Moody algebra. Let $\mathcal{E}\subset (\Hlie^*)^{\ZZ}$ be the subset of $c:\Hlie^*\rightarrow \ZZ$ such that $c(\lambda)=0$ for $\lambda$ outside the union of a finite number of sets of the form $D(\mu$). For $\lambda\in \Hlie^*$ denote $e(\lambda)\in\mathcal{E}$ such that $e(\lambda)(\mu)=\delta_{\lambda,\mu}$. $\mathcal{E}$ has a natural structure of commutative $\ZZ$-algebra such that $e(\lambda)e(\mu)=e(\lambda+\mu)$ (see \cite{kac}).

\noindent The classical character is the map $\text{ch}:\mathcal{O}(\U_q(\Glie))\rightarrow \mathcal{E}$ such that for $V\in\mathcal{O}(\U_q(\Glie))$:
$$\text{ch}(V)=\underset{\omega\in \Hlie^*}{\sum}\text{dim}(V_{\omega})e(\omega)$$
$\text{ch}$ is a ring morphism and $\text{ch}(L(\omega_1))=\text{ch}(L(\omega_2))\Rightarrow \omega_1=\omega_2$. 

\subsection{Formal character} Let $\U_q(\hat{\Glie})$ be a quantum affinization. In general the map $\text{ch}\circ\text{res}$ does not distinguish the simple integrable representations in $\mathcal{O}(\U_q(\hat{\Glie}))$. That is why Frenkel-Reshetikhin \cite{Fre} introduced the theory of $q$-characters for quantum affine algebras. We generalize the construction for quantum affinizations.

\noindent Let $V$ be in $\mathcal{O}(\U_q(\hat{\Glie}))$. For $\omega\in \Hlie^*$, the subspace $V_{\omega}\subset V$ is stable by the operators $\phi^{\pm}_{i,\pm m}$ ($i\in I$, $m\geq 0$). Moreover they commute and $[\phi_{i,m}^{\pm},k_h]=0$, so we have a pseudo-weight space decomposition: 
$$V_{\omega}=\underset{\gamma/ (\omega, \gamma)\in P_l}{\bigoplus} V_{\omega, \gamma}$$ 
where $V_{\omega, \gamma}$ is a simultaneous generalized eigenspace:
$$V_{\omega, \gamma}=\{x\in V_{\omega}/\exists p\in\NN,\forall i\in\{1,...,n\},\forall m\geq 0,(\phi_{i,\pm m}^{\pm}-\gamma_{i,\pm m}^{\pm})^p.x=0\}$$
As $V_{\omega}$ is finite dimensional the $V_{\omega, \gamma}$ are finite dimensional.

\noindent Let $\mathcal{E}_l\subset P_l^{\ZZ}$ be the ring of maps $c:P_l\rightarrow \ZZ$ such that $c(\lambda,\Psi)=0$ for $\lambda$ outside the union of a finite number of sets of the form $D(\mu$).

\begin{defi} The formal character of a module $V$ in the category $\mathcal{O}(\U_q(\hat{\Glie}))$ is $\text{ch}_q(V)\in \mathcal{E}_l$ defined by:
$$\text{ch}_q(V)=\underset{(\mu, \Gamma)\in P_l}{\sum} \text{dim}(V_{\mu, \Gamma}) e(\mu, \Gamma)$$
\end{defi}

\noindent We have the following commutative diagram:
$$\begin{array}{rcccl}
\mathcal{O}(\U_q(\hat{\Glie}))&\stackrel{\text{ch}_q}{\longrightarrow}&\mathcal{E}_l\\
\downarrow \text{res}&&\downarrow\beta\\
\mathcal{O}(\U_q(\Glie))&\stackrel{\text{ch}}{\longrightarrow}&\mathcal{E}\\\end{array}$$
where $\beta:\mathcal{E}_l\rightarrow\mathcal{E}$ is constructed from the first projection $\pi_1:P_l\rightarrow P$.

\subsection{Morphism of $q$-characters} The combinatorics of formal characters can be studied with a morphism of $q$-characters $\chi_q$ which is defined on a category $\mathcal{O}_{\text{int}}(\U_q(\hat{\Glie}))$ :

\subsubsection{The category $\mathcal{O}_{\text{int}}(\U_q(\hat{\Glie}))$}

Denote by $\mathcal{O}_{\text{int}}(\U_q(\Glie))$ (resp. $\mathcal{O}_{\text{int}}(\U_q(\hat{\Glie}))$) the category of integrable representations in the category $\mathcal{O}(\U_q(\Glie))$ (resp. $\mathcal{O}(\U_q(\hat{\Glie}))$). For example a simple integrable $l$-highest weight $\U_q(\hat{\Glie})$-modules is in $\mathcal{O}_{\text{int}}(\U_q(\hat{\Glie}))$. Moreover:

\begin{prop}\label{helpthm} For $V$ a module in $\mathcal{O}_{\text{int}}(\U_q(\hat{\Glie}))$ there are $P_{(\lambda,\Psi)}\geq 0$ ($(\lambda,\Psi)\in P_l^+$) such that:
$$\text{ch}_q(V)=\underset{(\lambda,\Psi)\in P_l^+}{\sum}P_{(\lambda,\Psi)}\text{ch}_q(L(\lambda,\Psi))$$
\end{prop}

\demo We have two preliminary points: 

1) a submodule, a quotient of an integrable module is integrable. 

2) for $V\in\mathcal{O}_{\text{int}}(\U_q(\hat{\Glie}))$ a module and $\mu$ a maximal weight of $V$, then there is $v\in V_{\mu}$ such that $\U_q(\hat{\Glie}).v$ is a $l$-highest weight module : indeed for $(\mu,\gamma)\in P_l$ such that $V_{\mu,\gamma}\neq \{0\}$ there is $v\in V_{\mu,\gamma}-\{0\}$ such that $\forall i\in I, r\geq 0$, $\phi_{i,\pm r}^{\pm}.v=\gamma_{i,\pm r}^{\pm} v$ (because for all $i\in I,r\geq 0$, $\text{Ker}(\phi_{i,\pm r}^{\pm}-\gamma_{i,\pm r}^{\pm})\cap V_{\mu,\gamma}\neq \{0\}$).

\noindent The end of the proof is essentially made in \cite{kac} (proposition 9.7) : first we prove that for $\lambda\in \Hlie^*$ there exists a filtration by a sequence of submodules in $\mathcal{O}_{\text{int}}(\U_q(\hat{\Glie}))$: $V=V_t\supset V_{t-1}\supset ... \supset V_1\supset V_0=0$ and $J\subset\{1,...,t\}$ such that:

(i) if $j\in J$, then $V_j/V_{j-1}\simeq L(\lambda_j,\Psi_j)$ for some $(\lambda_j,\Psi_j)\in P_l^+$ such that $\lambda_j\geq \lambda$

(ii) if $j\notin J$, then $(V_j/V_{j-1})_{\mu}=0$ for every $\mu\geq \lambda$

\noindent (see the lemma 9.6 of \cite{kac}). Next for $(\mu,\Psi)\in P_l^+$, fix $\lambda$ such that $\mu\geq \lambda$ and introduce $P_{(\mu,\Psi)}$ the number of times $(\mu,\Psi)$ appears among the $(\lambda_j,\Psi_j)$ (it is independent of the choice of the filtration and of $\mu$). We conclude as in proposition 9.7 of \cite{kac}.\qed

\begin{defi} $QP_l^+$ is the set of $(\mu,\gamma)\in P_l$ satisfying the following condition : 

i) there exist polynomials $Q_i(z), R_i(z)\in\CC[z]$ ($i\in I$) of constant term $1$ such that in $\CC[[z]]$ (resp. in $\CC[[z^{-1}]]$):
$$\underset{m\geq 0}{\sum}\gamma_{i,\pm m}^{\pm} z^{\pm m}= q_i^{\text{deg}(Q_i)-\text{deg}(R_i)}\frac{Q_i(zq_i^{-1})R_i(zq_i)}{Q_i(zq_i)R_i(zq_i^{-1})}$$

ii) there exist $\omega\in P^+$, $\alpha\in Q^+$ satisfying $\mu = \omega - \alpha$.\end{defi}

\noindent In particular $P_l^+\subset QP_l^+$.

\begin{prop}\label{fr} Let $V$ be a module in $\mathcal{O}_{\text{int}}(\U_q(\hat{\Glie}))$ and $(\mu,\gamma)\in P_l$. If $\text{dim}(V_{\mu, \gamma})>0$ then $(\mu,\gamma)\in QP_l^+$.\end{prop}

\demo The existence of the polynomials is shown as in \cite{Fre} (proposition 1): it reduces to the $sl_2$-case because for $v\in V$, $\hat{U}_i.v$ is finite dimensional. The existence of $\omega\in P$ and $\alpha\in Q^+$ is a consequence of proposition \ref{helpthm} and theorem \ref{simpint}.\qed

\subsubsection{Construction of $q$-characters} Consider formal variables $Y_{i,a}^{\pm}$ ($i\in I, a\in\CC^*$) and $k_{\omega}$ ($\omega\in\Hlie$). Let $\tilde{A}$ be the commutative group of monomials of the form $m=\underset{i\in I,a\in\CC^*}{\prod}Y_{i,a}^{u_{i,a}(m)} k_{\omega(m)}$ ($k_0=1$) where only a finite number of $u_{i,a}(m)\in\ZZ$ are non zero, $\omega(m)\in \Hlie$ (the coweight of $m$), and such that for $i\in I$:
$$\alpha_i(\omega(m))=r_iu_i(m)=r_i\underset{a\in\CC^*}{\sum}u_{i,a}(m)$$
The product is given by $u_{i,a}(m_1m_2)=u_{i,a}(m_1)+u_{i,a}(m_2)$ and $\omega(m_1m_2)=\omega(m_1)+\omega(m_2)$.

\noindent For example for $i\in I,a\in\CC^*$, we have $k_{\nu(\Lambda_i)}Y_{i,a}\in\tilde{A}$ because for $j\in I$, $\alpha_j(\nu({\Lambda_i}))=\Lambda_i(\nu(\alpha_j))=r_j\Lambda_i(\alpha_j^{\vee})=r_j\delta_{i,j}$. For $(\mu, \Gamma)\in QP_l^+$ we define $Y_{\mu,\Gamma}\in \tilde{A}$ by:
$$Y_{\mu, \Gamma}=k_{\nu(\mu)}\underset{i\in I, a\in\CC^*}{\prod}Y_{i,a}^{\beta_{i,a} -\gamma_{i,a}}$$
where $\beta_{i,a}, \gamma_{i,a}\in\ZZ$ are defined by $Q_i(u)=\underset{a\in\CC^*}{\prod}(1-ua)^{\beta_{i, a}}$ , $R_i(u)=\underset{a\in\CC^*}{\prod}(1-ua)^{\gamma_{i, a}}$. We have $Y_{\mu,\Gamma}\in\tilde{A}$ because for $i\in I$:
$$\alpha_i(\nu(\mu))=\mu(\nu(\alpha_i))=r_i\mu(\alpha_i^{\vee})=r_i(\text{deg}(Q_i)-\text{deg}(R_i))=r_iu_i(Y_{\mu,\Gamma})$$
For $\chi\in \tilde{A}^{\ZZ}$ we say $\chi\in\Yim$ if there is a finite number of element $\lambda_1,...,\lambda_s\in \Hlie^*$ such that the coweights of monomials of $\chi$ are in $\underset{j=1...s}{\bigcup}\nu(D(\lambda_j))$. In particular $\Yim$ has a structure of $\Hlie$-graded $\ZZ$-algebra.

\begin{defi} The $q$-character of a module $V\in\mathcal{O}_{\text{int}}(\U_q(\hat{\Glie}))$ is:
$$\chi_q(V)=\underset{(\mu,\Gamma)\in QP_l^+}{\sum} d(\mu,\Gamma)Y_{\mu,\Gamma}\in\Yim$$
where $d(\mu,\Gamma)\in\ZZ$ is defined by $\text{ch}_q(V)=\underset{(\mu,\Gamma)\in QP_l^+}{\sum} d(\mu,\Gamma)e(\mu,\Gamma)$.\end{defi}

\noindent We have a commutative diagram :
$$\begin{array}{rcccl}
\mathcal{O}_{\text{int}}(\U_q(\hat{\Glie}))&\stackrel{\chi_q}{\longrightarrow}&\Yim\\
\downarrow \text{res}&&\downarrow\beta\\
\mathcal{O}_{\text{int}}(\U_q(\Glie))&\stackrel{\text{ch}}{\longrightarrow}&\mathcal{E}\\\end{array}$$
where for $m\in\tilde{A}$, $\beta(m)=e(\omega(m))$.

\noindent If $C$ is of finite type then the weight of a monomial $m\in \Yim$ is $\omega(m)=\underset{i\in I}{\sum}u_i(m)\nu(\Lambda_i)$. So we can forget the $k_h$, and we get the $q$-characters defined in \cite{Fre}. In this case the integrable simple modules are finite dimensional.

\noindent Note that in the same way one can define the $q$-character of a finite dimensional $\U_q(\hat{\Hlie})$-module.

\subsubsection{Morphism of $q$-characters}\label{domdef}

\noindent Denote by $\text{Rep}(\U_q(\Glie))$ (resp. $\text{Rep}(\U_q(\hat{\Glie}))$) the Grothendieck group generated by the modules $V$ in $\mathcal{O}_{\text{int}}(\U_q(\Glie))$ (resp. $\mathcal{O}_{\text{int}}(\U_q(\hat{\Glie}))$) which have a composition series (a sequence of modules $V\supset V_1\supset V_2\supset ...$ such that $V_i/V_{i+1}$ is irreducible).

\noindent The tensor product defines a ring structure on $\text{Rep}(\U_q(\Glie))$ and $\text{ch}$ gives a ring morphism $\chi:\text{Rep}(\U_q(\Glie))\rightarrow \mathcal{E}$.

\noindent The $q$-characters are compatible with exact sequences and so we get a group morphism $\chi_q:\text{Rep}(\U_q(\hat{\Glie}))\rightarrow\Yim$ which is called morphism of $q$-characters.

\begin{prop} The morphism $\chi_q$ is injective and the following diagram is commutative:
$$\begin{array}{rcccl}
\text{Rep}(\U_q(\hat{\Glie}))&\stackrel{\chi_q}{\longrightarrow}&\Yim\\
\downarrow \text{res}&&\downarrow\beta\\
\text{Rep}(\U_q(\Glie))&\stackrel{\chi}{\longrightarrow}&\mathcal{E}\\\end{array}$$
\end{prop}

\noindent The commutativity of the diagram follows from the definition. To see that $\chi_q$ is injective, let us give some definitions:

\noindent A monomial $m\in\tilde{A}$ is said to be dominant if $u_{i,a}(m)\geq 0$ for all $i\in I,a\in\CC^*$. If a $l$-weight $(\omega,\Psi)$ belongs to $P_l^+$ then $Y_{(\omega,\Psi)}\in\tilde{A}$ is dominant. Moreover the map $(\omega,\Psi)\mapsto Y_{(\omega,\Psi)}$ defines a bijection between $P_l^+$ and dominant monomials. For $m\in\tilde{A}$ a dominant monomial we denote by $L(m)\in\Yim$ the $q$-character of $L(\omega,\Psi)$ where $(\omega,\Psi)$ is the corresponding dominant $l$-weight. In particular $L(m)=m+\text{monomials of lower weight}$ (in the sense of the ordering on $P$), and so the $L(m)$ are linearly independent.

\noindent A module with composition series in determined in the Grothendieck group by the multiplicity of the simple modules, and we have seen that the $\chi_q(L(\lambda,\Psi))$ ($(\lambda,\Psi)\in P_l^+$) are linearly independent in $\Yim$. So $\chi_q$ is injective.

\subsection{$q$-characters and universal $\mathcal{R}$-matrix} The original definition of $q$-characters (\cite{Fre}) was based on an explicit formula for the universal $\mathcal{R}$-matrix established in \cite{kt, lss, da}. In general no universal $\mathcal{R}$-matrix has been defined for a quantum affinization. However $q$-characters can be obtained with a piece of the formula of a ``$\mathcal{R}$-matrix'' in the same spirit as the original approach:

\noindent We refer to the chapter 3 of \cite{gui} for general background on $h$-formal deformations. Consider $\U_h(\hat{\Glie})$ the $\CC[[h]]$-algebra which is $h$-topologically generated by $\Hlie$ and the $x_{i,r}^{\pm}$ ($i\in I, r\in\ZZ$), $h_{i,m}$ ($i\in I, m\in\ZZ-\{0\}$) and with the relations of $\U_q(\hat{\Glie})$ (where we set for $\omega\in \Hlie$, $k_{\omega}=\text{exp}(h\omega)$). The subalgebra $\U_h(\hat{\Hlie})\subset \U_h(\hat{\Glie})$ is $h$-topologically generated by $\Hlie$ and the $h_{i,m}$ ($i\in I, m\in\ZZ-\{0\}$).

\noindent If $V$ is a $\U_q(\hat{\Glie})$-module (resp. $\U_q(\hat{\Hlie})$-module) which is $\U_q(\Hlie)$-diagonalizable then we have an algebra morphism $\pi_V(h):\U_h(\hat{\Glie})\rightarrow \text{End}(V)[[h]]$ (resp. $\pi_V(h):\U_h(\hat{\Hlie})\rightarrow \text{End}(V)[[h]]$) (Remark : for $\lambda\in\Hlie^*$, $\omega\in\Hlie$, $v\in V_{\lambda}$ we set $\omega.v=\lambda(\omega)v$).

\noindent Define $\mathcal{R}^0$ and $T$ in $\U_h(\hat{\Hlie})\hat{\otimes}\U_h(\hat{\Hlie})\subset \U_h(\hat{\Glie})\hat{\otimes}\U_h(\hat{\Glie})$ ($h$-topological completion of the tensor product) by the formula :
$$\mathcal{R}^0=\text{exp}(-(q-q^{-1})\underset{i,j\in I,m>0}{\sum}\frac{m}{[m]_q}\tilde{B}_{i,j}(q^m) h^m h_{i,m}\otimes h_{j,-m})$$
$$T=\text{exp}(-h\underset{1\leq i\leq 2n-l}{\sum}\Lambda_i^{\vee}\otimes \nu(\alpha_i))$$
Remark : we have the usual property of $T$ (see \cite{Fre}): for $\lambda,\mu\in\Hlie^*$, $x\in V_{\lambda}$, $y\in V_{\mu}$, we have $T.(x\otimes y)=q^{-(\lambda,\mu)}(x\otimes y)$. Indeed:
$$\underset{1\leq i\leq 2n-l}{\sum}\lambda(\Lambda_i^{\vee})\mu(\nu(\alpha_i))=(\mu,\underset{1\leq i\leq 2n-l}{\sum}\lambda(\Lambda_i^{\vee})\alpha_i)=(\mu,\lambda)$$
For $i\in I, m\in\ZZ-\{0\}$ denote $\tilde{h}_{i,m}=\underset{j\in I}{\sum}\tilde{C}_{j,i}(q^m)h_{j,m}$. We have an inclusion $\tilde{A}\subset \U_h(\hat{\Hlie})$ because the elements $Y_{i,a}^{\pm}=k_{\mp\nu(\Lambda_i)}\text{exp}(\mp(q-q^{-1})\underset{m\geq 1}{\sum}h^ma^{-m}\tilde{h}_{i,m})\in \U_h(\hat{\Glie})$ ($i\in I, a\in\CC^*$) are algebraically independent.

\begin{thm}\label{rmat} For $V$ a finite dimensional $\U_q(\hat{\Hlie})$-module, $((\text{Tr}_V\circ \pi_V(h))\otimes\text{Id})(\mathcal{R}^0T))\in \U_h(\hat{\Hlie})$ is equal to $\chi_q(V)$.\end{thm}

\demo For $(\lambda,\Psi)\in P_l$ consider $V_{(\lambda,\Psi)}$ and $((\text{Tr}_{V_{(\lambda,\Psi)}}\circ \pi_{V_{(\lambda,\Psi)}}(h))\otimes\text{Id})(\mathcal{R}^0T)$. First we see as in \cite{Fre} that the term $\mathcal{R}_0$ gives $\underset{i\in I, a\in\CC^*}{\prod}Y_{i,a}^{u_{i,a}(Y_{\lambda,\Psi})}$. But we have: 
$$\underset{1\leq i\leq 2n-l}{\sum}\lambda(\Lambda_i^{\vee})\nu(\alpha_i)=\nu(\underset{1\leq i\leq 2n-l}{\sum}\lambda(\Lambda_i^{\vee}) \alpha_i)=\nu(\lambda)$$
and so $T$ gives $k_{-\nu(\lambda)}$.\qed

\noindent In general for $V\in\mathcal{O}_{\text{int}}(\U_q(\hat{\Glie}))$ we can consider a filtration $(V_r)_{r\geq 0}$ of finite dimensional sub $\U_q(\hat{\Hlie})$-modules of $V$ such that $\underset{r\geq 0}{\bigcup}V_r=V$; so $\chi_q(V)$ is the ``limit'' of the $((\text{Tr}_{V_r}\circ \pi_{V_r}(h))\otimes\text{Id})(\mathcal{R}^0T)$ in $\Yim$.

\subsection{Combinatorics of $q$-characters}\label{comb} In this section we prove a symmetry property of general $q$-characters : the image of $\chi_q$ is the intersection of the kernels of screening operators (theorem \ref{sym}). Our proof is analog to the proof used by Frenkel-Mukhin \cite{Fre2} for quantum affine algebras; however new technical points are involved because of the $k_{\lambda}$ and infinite sums. In particular it shows that those $q$-characters are the combinatorial objects considered in \cite{her03} (which were constructed in the kernel of screening operators).

\noindent In sections \ref{comb} and \ref{dnc} we suppose that $C(z)$ is invertible (it includes the cases of quantum affine algebras and quantum toroidal algebras, see section \ref{bck}). We write $\tilde{C}(z)=\frac{\tilde{C}'(z)}{d(z)}$ where $d(z), \tilde{C}'_{i,j}(z)\in\ZZ[z^{\pm}]$. For $r\in\ZZ$ let $p_{i,j}(r)=[(D(z)\tilde{C}'(z))_{i,j}]_r$ where for a Laurent polynomial $P(z)\in\ZZ[z^{\pm}]$ we put $P(z)=\underset{r\in\ZZ}{\sum}[P(z)]_rz^r$.

\subsubsection{Construction of screening operators} Let $\Yim^{\text{int}}\subset \Yim$ be the subset consisting of those $\chi\in \Yim$ satisfying the following property : if $\lambda$ is the coweight of a monomial of $\chi$ there is $K\geq 0$ such that $k\geq K$ implies that for all $i\in I$, $\lambda \pm k r_i\alpha_i^{\vee}$ is not the coweight of a monomial of $\chi$.

\begin{lem} $\Yim^{\text{int}}$ is a subalgebra of $\Yim$ and $\text{Im}(\chi_q)\subset \Yim^{\text{int}}$.\end{lem}

\noindent Consider the free $\Yim^{\text{int}}$-module $\tilde{\Yim}_i=\underset{a\in\CC^*}{\prod}\Yim^{\text{int}}S_{i,a}$ and the linear map $\tilde{S}_i:\Yim^{\text{int}}\rightarrow \tilde{\Yim}_i$ such that, for a monomial $m$ :
$$\tilde{S}_i(m)=m\underset{a\in\CC^*}{\sum}u_{i,a}(m)S_{i,a}$$
In particular $\tilde{S}_i$ is a derivation. Let us choose a representative $a$ for each class of $\CC^*/q_i^{2\ZZ}$ and consider:
$$\Yim_i=\underset{a\in\CC^*/q_i^{2 \ZZ}}{\prod}\Yim^{\text{int}}S_{i,a}$$
For $i\in I$ and $a\in\CC^*$ we set:
$$A_{i,a}=k_iY_{i,aq_i^{-1}}Y_{i,aq_i}\underset{j/C_{j,i}<0\text{ , }r=C_{j,i}+1,C_{j,i}+3,...,-C_{j,i}-1}{\prod}Y_{j,aq^r}^{-1}\in\tilde{A}$$
We have $A_{i,a}\in\tilde{A}$ because for $j\in I$: $\alpha_j(r_i\alpha_i^{\vee})=r_iC_{i,j}=r_jC_{j,i}=r_ju_j(A_{i,a})$.

We would like to see $\Yim_i$ as a quotient of $\tilde{\Yim}_i$ by the relations $S_{i,aq_i}=A_{i,a}S_{i,aq_i^{-1}}$. But the projection is not defined for all elements of $\tilde{\Yim}_i$ because there are infinite sums. However if $\chi\in \Yim^{\text{int}}$ and $m$ is a monomial of $\chi$ there is a finite number of monomials in $\chi$ of the form $mA_{i,aq_i}^{-1}A_{i,aq_i^3}^{-1}...A_{i,aq_i^r}^{-1}$ or of the form $mA_{i,aq_i^{-1}}A_{i,aq_i^{-3}}^{-1}...A_{i,aq_i^{-r}}^{-1}$. So the projection on $\Yim_i$ is well defined on $\tilde{S}_i(\Yim^{\text{int}})\subset \tilde{\Yim}_i$. In particular we can define by projection of $\tilde{S}_i$ the $i^{th}$ screening operator $S_i:\Yim^{\text{int}}\rightarrow \Yim_i$.

\noindent The original definition for the finite case is in \cite{Fre}.

\subsubsection{The morphism $\tau_i$} Some operators $\tau_i$ ($i\in I$) were defined for the finite case in \cite{Fre2}. We generalize the construction and the properties of the operators $\tau_i$ (lemma \ref{aideun} and \ref{aidedeux}).

\noindent Let $i\in I$. Denote $\Hlie_i^{\perp}=\{\omega\in \Hlie/\alpha_i(\omega)=0\}$. 

\noindent Consider formal variables $k^{(i)}_{r}$ ($r\in\ZZ$), $k_{\omega}$ ($\omega\in\Hlie$), $Y_{i,a}^{\pm}$ ($a\in\CC^*$), $Z_{j,c}$ ($j\in I-\{i\}$, $c\in\CC^*$). Let $\tilde{A}^{(i)}$ be the commutative group of monomials :
$$m=k^{(i)}_{r(m)}k_{\omega(m)}\underset{a\in\CC^*}{\prod}Y_{i,a}^{u_{i,a}(m)}\underset{j\in I,j\neq i,c\in\CC^*}{\prod}Z_{j,c}^{z_{j,c}(m)}$$ 
where only a finite number of $u_{i,a}(m),z_{j,c}(m),r(m)\in\ZZ$ are non zero, $\omega(m)\in \Hlie_i^{\perp}$ and such that $r(m)=r_iu_i(m)=r_i\underset{a\in\CC^*}{\sum}u_{i,a}(m)$. The product is defined as for $\tilde{A}$. We call $(r(m),\omega(m))\in\ZZ\times \Hlie_i^{\perp}$ the coweight of the monomial $m$.

\noindent Let $\tau_i:\tilde{A}\rightarrow \tilde{A}^{(i)}$ be the group morphism defined by ($j\in I$, $a\in\CC^*$, $\lambda\in \Hlie$):
$$\tau_i(Y_{j,a})=Y_{i,a}^{\delta_{i,j}}\underset{k\neq i\text{ , }r\in\ZZ}{\prod}Z_{k,a q^r}^{p_{j,k}(r)}\text{ , }\tau_i(k_{\lambda})=k^{(i)}_{\alpha_i(\lambda)}k_{\lambda-\alpha_i(\lambda)\frac{\alpha_i^{\vee}}{2}}$$
(note that it is a formal definition because $Y_{j,a}k_{\nu(\Lambda_j)}\in\tilde{A}$ but $Y_{j,a}\notin \tilde{A}$). It is well defined because for $m\in \tilde{A}$, $\alpha_i(\omega(m))=r_iu_i(m)$ and $\alpha_i(\omega(m)-\alpha_i(\omega(m))\frac{\alpha_i^{\vee}}{2})=0$.

\begin{lem}\label{aideun} The morphism $\tau_i$ is injective and for $a\in\CC^*$ we have:
$$\tau_i(A_{i,a})=k_{2r_i}^{(i)}Y_{i,aq_i^{-1}}Y_{i,aq_i}$$\end{lem}

\demo Let $m\in\tilde{A}$ such that $\tau_i(m)=1$. For $a\in\CC^*$ we have $u_{i,a}(m)=u_{i,a}(\tau_i(m))=0$. For $k\in I$, $a\in\CC^*$ denote $u_{k,a}(m)(z)=\underset{r\in\ZZ}{\sum}u_{k,aq^r}(m)z^r\in\ZZ[z^{\pm}]$. For $j\in I - \{i\}$, we have :
$$0=z_{j,aq^R}(\tau_i(m))=\underset{k\in I,r+r'=R}{\sum}p_{k,j}(r')u_{k,aq^r}(m)=[\underset{k\in I}{\sum}\tilde{C}'_{k,j}(z)u_{k,a}(m)(z)]_R$$
As $\tilde{C}(z)$ is invertible we get $u_{k,a}(m)=0$ for all $a\in\CC^*$. In particular for $j\in I$ we have $\alpha_j(\omega(m))=r_ju_j(m)=0$. But $\omega(m)-\alpha_i(\omega(m))\frac{\alpha_i^{\vee}}{2}=0=\omega(m)$ and so $m=1$.

\noindent For the second point let $M=\tau_i(A_{i,a})$. First for $b\in\CC^*$, $u_{i,b}(M)=u_{i,b}(A_{i,a})=\delta_{a/b,q_i}+\delta_{a/b,q_i^{-1}}$. For $R\in\ZZ$ and $j\neq i$ we have:
$$z_{j,aq^R}(M)=[(\tilde{C}'(z)C(z))_{i,j}]_R=[(d(z)D(z))_{i,j}]_R=0$$
Finally we have $r(M)=r_i\alpha_i(\alpha_i^{\vee})=-2r_i$ and $\omega(M)=r_i\alpha_i^{\vee}-r_i\alpha_i(\alpha_i^{\vee})\frac{\alpha_i^{\vee}}{2}=0$.\qed

\noindent Formally we have $\tau_i(k_i)=k_{2r_i}^{(i)}$ and for $j\in I-\{i\}$ $\tau_i(k_j)=k_{B_{j,i}}^{(i)}k_{\alpha_j^{(i)}}$ where $\alpha_j^{(i)}=r_j\alpha_j^{\vee}-\frac{B_{j,i}}{2}\alpha_i^{\vee}$. This motivates the following definition: for $(r,\omega)\in \ZZ\times \Hlie_i^{\perp}$ denote :
$$D(r,\omega)=\{(r',\omega')\in\ZZ\times \Hlie_i^{\perp}/\omega'=\omega-\underset{j\in I,j\neq i}{\sum}m_j\alpha_j^{(i)}\text{ , }r'=r-\underset{j\in I, j\neq i}{\sum}B_{j,i}m_j-2r_ik/m_j,k\geq 0\}$$

\noindent Define $\Yim^{\text{int}, (i)}\subset (\tilde{A}^{(i)})^{\ZZ}$ as the set of $\chi$ such that : 

i) there is a finite number of elements $(r_1,\omega_1),...,(r_s,\omega_s)\in \ZZ\times \Hlie_i^{\perp}$ such that the coweights of monomials of $\chi$ are in $\underset{j=1...s}{\bigcup}D(r_j,\omega_j)$.

ii) for $(r,\lambda)$ the coweight of a monomial of $\chi$ there is $K\geq 0$ such that $k\geq K$ implies that for all $j\in I$, $j\neq i$, $(r\pm B_{j,i}k,\lambda \pm k\alpha_j^{(i)})$ and $(r\pm 2kr_i,\lambda)$ are not the coweight of a monomial of $\chi$.

\noindent In particular $\Yim^{\text{int}, (i)}$ has a structure of $\ZZ\times \Hlie_i^{\perp}$-graded $\ZZ$-algebra.

\noindent The morphism $\tau_i$ can be extended to a unique morphism of $\ZZ$-algebra $\tau_i:\Yim^{\text{int}}\rightarrow \Yim^{\text{int}, (i)}$. Denote by $\chi_q^i$ the morphism of $q$-characters for the algebra $\U_{q_i}(\hat{sl_2})$.

\begin{lem}\label{aidedeux} Consider $V\in\mathcal{O}_{\text{int}}(\U_q(\hat{\Glie}))$ and a decomposition $\tau_i(\chi_q(V))=\underset{k}{\sum}P_kQ_k$ where 
\\$P_k\in\ZZ[Y_{i,a}^{\pm}k_{\pm r_i}^{(i)}]_{a\in\CC^*}$, $Q_k$ is a monomial in $\ZZ[Z_{j,c}^{\pm}, k_h]_{j\neq i,a\in\CC^*,h\in \Hlie_i^{\perp}}$ and all monomials $Q_k$ are distinct. Then there exists a $\hat{U}_i$-module $\underset{k}{\bigoplus}V_k$ isomorphic to the restriction of $V$ to $\hat{U}_i$ and such that $\chi_q^i(V_k)=P_k$.\end{lem}

\demo Let $\U_q(\hat{\Hlie})_i^{\perp}$ the subalgebra of $\U_q(\hat{\Glie})$ generated by the $k_h$ ($h\in \Hlie_i^{\perp}$), $h_{j,m}$ ($j\neq i$, $m\in\ZZ-\{0\}$). We can apply the proof of lemma 3.4 of \cite{Fre2} with $\hat{U}_i$ and $\U_q(\hat{\Hlie})_i^{\perp}$ because :

i) $\hat{U}_i$ and $\U_q(\hat{\Hlie})_i^{\perp}$ commute in $\U_q(\hat{\Glie})$ 

ii) The image $\omega-\alpha_i(\omega)\frac{\alpha_i^{\vee}}{2}$ in $\Hlie_i^{\perp}$ of $\omega\in \Hlie$ suffices to encode the action of the $k_h$ ($h\in \Hlie_i^{\perp}$) on a vector of weight $\nu^{-1}(\omega)=\lambda$. Indeed for $h\in \Hlie_i^{\perp}$, we have:
$$\lambda(h)=(\nu^{-1}(h),\nu^{-1}(\omega))=\nu^{-1}(h)(\omega)=\nu^{-1}(h)(\omega-\omega(\alpha_i)\frac{\alpha_i^{\vee}}{2})$$
because $\alpha_i(h)=0\Rightarrow \nu^{-1}(\alpha_i^{\vee})=0$.\qed

\subsubsection{$\tau_i$ and screening operators} In this section we prove that $\text{Im}(\chi_q)\subset\text{Ker}(S_i)$ (proposition \ref{but}) with a generalization of the proof of Frenkel-Mukhin \cite{Fre2}.

\noindent Consider the $\Yim^{\text{int}, (i)}$-module $\tilde{\Yim}^{(i)}_i=\underset{a\in\CC^*}{\prod}\Yim^{\text{int}, (i)}S_{i,a}$ and the linear map $\overline{S}_i:\Yim^{\text{int}, (i)}\rightarrow \tilde{\Yim}_i^{(i)}$ such that, for a monomial $m$ :
$$\overline{S}_i(m)=m\underset{a\in\CC^*}{\sum}u_{i,a}(m)S_{i,a}$$
In particular $\overline{S}_i$ is a derivation. Consider $\Yim_i^{(i)}=\underset{a\in\CC^*/q_i^{2 \ZZ}}{\prod}\Yim^{\text{int}, (i)}S_{i,a}$. The $\overline{S}_i(\Yim^{\text{int}, (i)})\subset \tilde{\Yim}_i^{(i)}$ can be projected in $\Yim_i^{(i)}$ by the relations :
$$S_{i,aq_i}=Y_{i,aq_i}Y_{i,aq_i^{-1}}k_{2r_i}^{(i)}S_{i,aq_i^{-1}}$$
and we get a derivation that we denote also by $\overline{S}_i:\Yim^{\text{int}, (i)}\rightarrow \Yim_i^{(i)}$.

\noindent We also define a map $\tau_i:\Yim_i\rightarrow \Yim_i^{(i)}$ in an obvious way (with the help of lemma \ref{aideun}). We see as in lemma 5.4 of \cite{Fre2} that:
\begin{lem}\label{comm} We have a commutative diagram: 
$$\begin{array}{rcccl}
\Yim^{\text{int}}&\stackrel{S_i}{\longrightarrow}&\Yim_i\\
\downarrow \tau_i&&\downarrow\tau_i\\
\Yim^{\text{int},(i)}&\stackrel{\overline{S}_i}{\longrightarrow}&\Yim_i^{(i)}\\\end{array}$$
\end{lem}

\noindent With the help of lemma \ref{aideun}, \ref{aidedeux} and \ref{comm} we see as in corollary 5.5 of \cite{Fre2}:

\begin{lem}\label{but} We have $\text{Im}(\chi_q)\subset\underset{i\in I}{\bigcap}\text{Ker}(S_i)$.\end{lem}

\noindent In the following we denote $\mathfrak{K}_i=\text{Ker}(S_i)$ and $\mathfrak{K}=\underset{i\in I}{\bigcap}\mathfrak{K}_i$.

\begin{lem}\label{trois} An element $\chi\in\Yim^{\text{int}}$ is in $\mathfrak{K}_i$ if and only if it can be written in the form $\chi=\underset{k}{\sum} P_kQ_k$ where $P_k\in \ZZ[k_{\nu(\Lambda_i)}Y_{i,a}(1+A_{i,aq_i}^{-1})]_{a\in\CC^*}$, $Q_k$ is a monomial in $\ZZ[Y_{j,a}^{\pm},k_h]_{j\neq i,a\in\CC^*,h\in P_i^{*,\perp}}$ and all monomials $Q_k$ are distinct.\end{lem}

\demo We use the result for the $sl_2$-case which is proved in \cite{Fre}. First an element of this form is in $\mathfrak{K}_i$. Consider $\chi\in\mathfrak{K}_i$ and write $\tau_i(\chi)=\underset{k}{\sum} P_k'Q_k'$ as in lemma \ref{aidedeux}. From lemma \ref{comm} we have $0=\overline{S}_i(\chi)=\underset{k}{\sum}\overline{S}_i(P_k')Q_k$. So all $\overline{S}_i(P_k')=0$ and it follows from the $sl_2$-case that $P_k'\in \ZZ[Y_{i,a}k_{r_i}^{(i)}+Y_{i,aq_i^2}^{-1}k_{-r_i}^{(i)}]_{a\in\CC^*}$. The lemma \ref{aideun} lead us to the conclusion.\qed

\subsubsection{Description of $\text{Im}(\chi_q)$}

Dominant monomials are defined in section \ref{domdef}. We have:

\begin{lem}\label{uni} An element $\chi\in\mathfrak{K}$ has at least one dominant monomial.\end{lem}

\noindent With the help of lemma \ref{trois} we can use the proof of lemma 5.6 of \cite{Fre2} (see also the proof of theorem 4.9 in \cite{her01} at $t=1$).

\begin{thm}\label{sym} We have $\text{Im}(\chi_q)=\mathfrak{K}$. Moreover the elements of $\mathfrak{K}$ are the sums:
$$\underset{m\text{ dominant}}{\sum}\lambda_m L(m)$$
where $\lambda_m=0$ for $\omega(m)$ outside the union of a finite number of sets of the form $D(\mu$).\end{thm}

\demo The inclusion $\text{Im}(\chi_q)\subset\mathfrak{K}$ is proved in lemma \ref{but}. For the other one, consider $\chi\in\mathfrak{K}$. We can suppose that the weights of $\chi$ are in a set $D(\lambda)$ (because the weights of each $L(m)$ are in a set $D(\mu)$). We define by induction $L^{(k)}(m)\in\text{Im}(\chi_q)$ ($k\geq 0$) in the following way: we set $L^{(0)}=\underset{\omega(m)=\lambda}{\sum}[\chi]_m L(m)$. If $L^{(k)}$ is defined, we consider the set $\tilde{A}_{k+1}$ of monomials $m'$ which appear in $\chi-L^{(k)}$ such that $\lambda-\omega(m')=m_1r_1\alpha_1^{\vee}+...+m_nr_n\alpha_n^{\vee}$ where $m_1,...,m_n\geq 0$ and $m_1+...+m_n=k$. We set:
$$L^{(k+1)}=L^{(k)}+\underset{m'\in \tilde{A}_{k+1}}{\sum}[\chi-L^{(k)}]_{m'}L(m')$$
Then we set $L^{\infty}=\underset{k\geq 0/m\in \tilde{A}_{k}}{\sum}[L^{(k)}]_mL(m)\in\text{Im}(\chi_q)$ and it follows from lemma \ref{uni} that $L^{\infty}=\chi$.\qed

\noindent Note that proposition \ref{helpthm} gives that for $\chi_q(V)$ ($V$ module in $\mathcal{O}_{\text{int}}(\U_q(\hat{\Glie}))$) the $\lambda_m$ are non-negative.

\noindent Remark: for $m\in \tilde{A}$ a dominant monomial we prove in the same way that there is a unique $F(m)\in\mathfrak{K}$ such that $m$ has coefficient $1$ in $F(m)$ and $m$ is the unique dominant monomial in $F(m)$. In the finite case an algorithm was given by Frenkel-Mukhin \cite{Fre2} to compute the $F(m)$. In \cite{her03} we extended the definition of the algorithm for generalized Cartan matrix and showed that it is well-defined if $i\neq j\Rightarrow C_{i,j}C_{j,i}\leq 3$ (see also \cite{her02} for the detailed description of this algorithm at $t=1$). Theorem \ref{sym} allows us to prove two results announced in \cite{her03} : the algorithm is well defined for 

$A_1^{(1)}$ (with $r_1=r_2=2$) because $\text{det}(C(z))=z^4-z^2-z^{-2}+z^{-4}\neq 0$

$A_2^{(2)}$ (with $r_1=4$, $r_2=1$) because $\text{det}(C(z))=z^5-z-z^{-1}+z^{-5}\neq 0$

\noindent But for $A_1^{(1)}$ (with $r_1=r_2=1$) we have $\text{det}(C(z))=0$; we observed in \cite{her03} that the algorithm is not well-defined in this case.

\section{Drinfel'd new coproduct and fusion product}\label{dnc} Our study of combinatorics of $q$-characters gives a ring structure on $\text{Im}(\chi_q)$ (corollary \ref{ring}). As $\chi_q$ is injective we get an induced ring structure on the Grothendieck group. In this section we prove that it is a fusion product (theorem \ref{posc}), that is to say that the product of two modules is a module. We use the Drinfel'd new coproduct (proposition \ref{coprod}); as it involves infinite sums, we have to work in a larger category where the tensor product is well-defined (theorem \ref{prod}). To end the proof of theorem \ref{posc} we define specializations of certain forms which allow us to go from the larger category to $\mathcal{O}(\U_q(\hat{\Glie}))$. Note that in our construction we do not use that $C(z)$ is invertible.

\subsection{Fusion product}

As the $S_i$ are derivations, theorem \ref{sym} gives:

\begin{cor}\label{ring} $\text{Im}(\chi_q)$ is a subring of $\Yim$.\end{cor}

\noindent Since $\chi_q$ is injective on $\text{Rep}(\U_q(\hat{\Glie}))$, the product of $\Yim$ gives an induced commutative product $*$ on $\text{Rep}(\U_q(\hat{\Glie}))$. For $(\lambda,\Psi), (\lambda',\Psi')\in P_l^+$ there are $Q_{\lambda,\Psi,\lambda',\Psi'}(\mu,\Phi)\in\ZZ$ such that:
$$L(\lambda,\Psi)*L(\lambda',\Psi')=L(\lambda+\lambda',\Psi\Psi')+\underset{(\mu,\Phi)\in P_l^+/ \mu < \lambda +\lambda'}{\sum}Q_{\lambda,\Psi,\lambda',\Psi'}(\mu,\Phi)L(\mu,\Phi)$$
We will interpret this product as a fusion product related to the basis of simple modules : that is to say we will show that a product of modules is a module (see \cite{Fu} for generalities on fusion rings and physical motivations). Let us explain it in more details : consider :
$$\text{Rep}^+(\U_q(\hat{\Glie}))=\underset{(\lambda,\Psi)\in P_l^+}{\bigoplus}\NN.L(\lambda,\Psi)\subset\text{Rep}(\U_q(\hat{\Glie}))=\underset{(\lambda,\Psi)\in P_l^+}{\bigoplus}\ZZ.L(\lambda,\Psi)$$

\begin{thm}\label{posc} The subset $\text{Rep}^+(\U_q(\hat{\Glie}))\subset\text{Rep}(\U_q(\hat{\Glie}))$ is stable by *.\end{thm}

\noindent In this section \ref{dnc} we prove this theorem by interpreting $*$ with the help of a generalization of the new Drinfel'd coproduct. Note that theorem \ref{posc} means that for $(\lambda,\Psi), (\lambda',\Psi')\in P_l^+$ we have $Q_{\lambda,\Psi,\lambda',\Psi'}(\mu,\Phi)\geq 0$.

\subsection{Coproduct} 

\subsubsection{Reminder: case of a quantum affine algebra and Drinfel'd-Jimbo coproduct}

As said before the case of a quantum affine algebra is a very special one because there are two realizations (if we add a central charge); in particular there is a coproduct on $\U_q(\hat{\Glie})$, a tensor product on $\mathcal{O}_{\text{int}}(\U_q(\hat{\Glie}))$ and $\text{Rep}(\U_q(\hat{\Glie}))$ is a ring. It is the product * because it is shown in \cite{Fre} that $\chi_q$ is a ring morphism. In particular the tensor product is commutative. So theorem \ref{posc} is proved in this case.

\subsubsection{General case : new Drinfel'd coproduct} In general we have a coproduct $\Delta_{\hat{\Hlie}}:\U_q(\hat{\Hlie})\rightarrow \U_q(\hat{\Hlie})\otimes\U_q(\hat{\Hlie})$ for the commutative algebra $\U_q(\hat{\Hlie})$ defined by ($h\in P^*$, $i\in I$, $m\neq 0$):
$$\Delta_{\hat{\Hlie}}(k_h)=k_h\otimes k_h\text{ , }\Delta_{\hat{\Hlie}}(h_{i,m})=1\otimes h_{i,m} + h_{i,m}\otimes 1$$
In particular we have ($i\in I, m\geq 0$) : $\Delta_{\hat{\Hlie}}(\phi_{i,\pm m}^{\pm})= \underset{0\leq l\leq m}{\sum}\phi_{i,\pm (m-l)}^{\pm}\otimes\phi_{i,\pm l}^{\pm}$.

\noindent No coproduct has been defined for the entire $\U_q(\hat{\Glie})$. However Drinfel'd (unpublished note, see also \cite{di, df}) defined for $\U_q(\hat{sl_n})$ a map which behaves as a new coproduct adapted to the affinization realization. In this section we use those formulas for general quantum affinizations; as infinite sums are involved, we use a formal parameter $u$ so that it makes sense.

\noindent Let $\mathcal{C}=\CC((u))$ be the field of Laurent series $\underset{r\geq R}{\sum}\lambda _r u^r$ ($R\in\ZZ$, $\lambda_r\in\CC$). The algebra $\tilde{\U}_q(\hat{\Glie})$ is defined in section \ref{stat}. Consider the $\mathcal{C}$-algebra $\tilde{\U}_q'(\hat{\Glie})=\mathcal{C}\otimes\tilde{\U}_q(\hat{\Glie})$ (resp. $\U_q'(\hat{\Glie})=\mathcal{C}\otimes\U_q(\hat{\Glie})$). Let $\tilde{\U}_q'(\hat{\Glie})\hat{\otimes}\tilde{\U}_q'(\hat{\Glie})=(\tilde{\U}_q(\hat{\Glie})\otimes_{\CC}\tilde{\U}_q(\hat{\Glie}))((u))$ be the $u$-topological completion of $\tilde{\U}_q'(\hat{\Glie})\otimes_{\mathcal{C}}\tilde{\U}_q'(\hat{\Glie})$. It is also a $\mathcal{C}$-algebra.

\begin{prop}\label{coprod} There is a unique morphism of $\mathcal{C}$-algebra $\Delta_u: \tilde{\U}_q'(\hat{\Glie})\rightarrow \tilde{\U}_q'(\hat{\Glie})\hat{\otimes}\tilde{\U}_q'(\hat{\Glie})$ such that for $i\in I$, $r\in\ZZ$, $m\geq 0$, $h\in\Hlie$:
$$\Delta_u(x_{i,r}^+)= x_{i,r}^+\otimes 1 + \underset{l\geq 0}{\sum} u^{r+l} (\phi_{i,-l}^-\otimes x_{i,r+l}^+)$$
$$\Delta_u(x_{i,r}^-)=u^r (1\otimes x_{i,r}^-) + \underset{l\geq 0}{\sum} u^l (x_{i,r-l}^-\otimes \phi_{i,l}^+)$$
$$\Delta_u(\phi_{i,\pm m}^{\pm})= \underset{0\leq l\leq m}{\sum}u^{\pm l} (\phi_{i,\pm (m-l)}^{\pm}\otimes\phi_{i,\pm l}^{\pm})\text{ , }\Delta_u(k_h)=k_h\otimes k_h$$
\end{prop}

\demo We can easily check the compatibility with relations (\ref{afcartu}), (\ref{actcartc}), (\ref{actcartplus}), (\ref{actcartmoins}), (\ref{partact}), (\ref{plusmoinsc}) because $\Delta_u$ can also be given in terms of the currents of section \ref{cur}: we have in $(\tilde{\U}_q'(\hat{\Glie})\otimes_{\mathcal{C}}\tilde{\U}_q'(\hat{\Glie}))[[z,z^{-1}]]$ :
$$\Delta_u(x_i^+(z))=x_i^+(z)\otimes 1 + \phi_i^-(z)\otimes x_i^+(zu)\text{ , }\Delta_u(x_i^-(z))=1\otimes x_i^-(zu) + x_i^-(z)\otimes \phi_i^+(zu)$$
$$\Delta_u(\phi_i^{\pm}(z))=\phi_i^{\pm}(z)\otimes\phi_i^{\pm}(zu)$$\qed

\noindent Remark 1 : If $C$ is finite or simply laced then $\Delta_u$ is compatible with the affine quantum Serre relations (relations (\ref{equadeux})) and can be defined for $\U_q'(\hat{\Glie})$ (see \cite{di} for finite symmetric cases and \cite{e, gr} for other finite cases). We conjecture that it is also true for general $C$, but we do not need it for our purposes.

\noindent Remark 2 : let $T:\tilde{\U}_q(\hat{\Glie})\rightarrow\tilde{\U}_q'(\hat{\Glie})$ be the $\ZZ$-gradation morphism defined by $T(x_{i,r}^{\pm})=u^rx_{i,r}^{\pm}$, $T(\phi_{i,m}^{\pm})=u^m\phi_{i,m}^{\pm}$, $T(k_h)=k_h$. The $u$ is put in such a way that $\Delta_u=(\text{Id}\otimes T)\circ\Delta$ where $\Delta$ is the usual new Drinfel'd coproduct (without $u$).

\noindent Remark 3 : The map $\Delta_u$ is not coassociative , indeed in $(\tilde{\U}_q'(\hat{\Glie})\otimes_{\mathcal{C}}\tilde{\U}_q'(\hat{\Glie})\otimes_{\mathcal{C}}\tilde{\U}_q'(\hat{\Glie}))[z]$:
$$((\Delta_u\otimes \text{Id})\circ \Delta_u)(\phi_i^+(z))=\phi_i^+(z)\otimes \phi_i^+(uz)\otimes \phi_i^+(uz)$$
$$((\text{Id}\otimes \Delta_u)\circ \Delta_u)(\phi_i^+(z))=\phi_i^+(z)\otimes \phi_i^+(uz)\otimes \phi_i^+(u^2z)$$

\noindent Remark 4 : Although is is not defined in a strict sense, the ``limit'' of $\Delta_u$ at $u=1$ is coassociative. On $\U_q(\hat{\Hlie})$ the limit at $u=1$ makes sense and is $\Delta_{\hat{\Hlie}}$.

\subsection{Tensor products of representations of $\tilde{\U}_q'(\hat{\Glie})$} As the coproduct involves infinite sums, we have to introduce a category larger than $\mathcal{O}(\U_q(\hat{\Glie}))$ in order to define tensor products:

\subsubsection{The category $\mathcal{O}(\tilde{\U}_q'(\hat{\Glie}))$}

\begin{defi} The set of $l,u$-weights $P_{l,u}$ is the set of couple $(\lambda,\Psi(u))$ such that $\lambda\in \Hlie^*$, $\Psi(u)=(\Psi_{i,\pm m}^{\pm}(u))_{i\in I, m\geq 0}$, $\Psi_{i,\pm m}^{\pm}(u)\in\CC[u,u^{-1}]$ and $\Psi_{i,0}^{\pm}(u)=q_i^{\pm \lambda(\alpha_i^{\vee})}$.\end{defi}

\begin{defi}\label{oprime} An object $V$ of the category $\mathcal{O}(\tilde{\U}_q'(\hat{\Glie}))$ is a $\mathcal{C}$-vector space with a structure of $\tilde{\U}_q'(\hat{\Glie})$-module such that:

i) $V$ is $\U_q(\Hlie)$-diagonalizable

ii) For all $\lambda\in\Hlie^*$, the sub $\mathcal{C}$-vector space $V_{\lambda}\subset V$ is finite dimensional

iii) there are a finite number of element $\lambda_1,...,\lambda_s\in \Hlie^*$ such that the weights of $V$ are in $\underset{j=1...s}{\bigcup}D(\lambda_j)$

iv) for $\lambda\in\Hlie^*$, $V_{\lambda}=\underset{(\lambda,\Psi(u))\in P_{l,u}}{\bigoplus}V_{(\lambda,\Psi(u))}$ where:
$$V_{\lambda, \Psi(u)}=\{x\in V_{\lambda}/\exists p\in\NN,\forall i\in\{1,...,n\},\forall r\geq 0,(\phi_{i,\pm r}^{\pm}-\Psi_{i,\pm r}^{\pm}(u))^p.x=0\}$$\end{defi}

\noindent The property iv) is added because $\mathcal{C}$ is not algebraically closed.

\noindent The scalar extension and the projection $\tilde{\U}_q(\hat{\Glie})\rightarrow \U_q(\hat{\Glie})$ gives an injection $i:\mathcal{O}(\U_q(\hat{\Glie}))\rightarrow \mathcal{O}(\tilde{\U}_q'(\hat{\Glie}))$ such that for $V\in \mathcal{O}(\U_q(\hat{\Glie}))$, $i(V)=V\otimes\mathcal{C}$.

\noindent Let $\mathcal{E}_{l,u}\subset P_{l,u}^{\ZZ}$ be defined as $\mathcal{E}_l$. The formal character of a module $V$ in the category $\mathcal{O}(\tilde{\U}_q'(\hat{\Glie}))$ is:
$$\text{ch}_{q,u}(V)=\underset{(\mu, \Gamma(u))\in P_{l,u}}{\sum} \text{dim}_{\mathcal{C}}(V_{\mu, \Gamma(u)}) e(\mu, \Gamma(u))\in\mathcal{E}_{l,u}$$
We have a map $i_{\mathcal{E}}:\mathcal{E}_l\rightarrow\mathcal{E}_{l,u}$ such that $i_{\mathcal{E}}((\lambda,\Psi))=(\lambda,(\Psi_{i,\pm m}^{\pm m}))$ and a commutative diagram:
$$\begin{array}{rcccl}
\mathcal{O}(\U_q(\hat{\Glie}))&\stackrel{\text{ch}_q}{\longrightarrow}&\mathcal{E}_l\\
\downarrow \text{i}&&\downarrow i_{\mathcal{E}}\\
\mathcal{O}(\tilde{\U}_q'(\hat{\Glie}))&\stackrel{\text{ch}_{q,u}}{\longrightarrow}&\mathcal{E}_{l,u}\\\end{array}$$
In an analogous way one defines the category $\mathcal{O}(\U_q'(\hat{\Glie}))$ and a formal character $\text{ch}_{q,u}$ on $\mathcal{O}(\U_q'(\hat{\Glie}))$.

\subsubsection{Tensor products} We consider subcategories of $\mathcal{O}(\tilde{\U}_q'(\hat{\Glie}))$. Let $R\in\ZZ$, $R\geq 0$:

\begin{defi} $\mathcal{O}^R(\tilde{\U}_q'(\hat{\Glie}))$\label{orfus} is the category of modules $V\in\mathcal{O}(\tilde{\U}_q'(\hat{\Glie}))$ such that for all $\lambda\in\Hlie^*$, there is a $\mathcal{C}$-basis $(v_{\alpha}^{\lambda})_{\alpha}$ of $V_{\lambda}$ satisfying :

i) for all $m\in\ZZ$, $\alpha$, $\beta$, the coefficient of $x_{i,m}^+.v_{\alpha}^{\lambda}$ on $v_{\beta}^{\lambda+\alpha_i}$ (resp. of $x_{i,m}^-.v_{\alpha}^{\lambda}$ on $v_{\beta}^{\lambda-\alpha_i}$) is in $\CC[[u]]$ if $m\geq 0$, in $u^{Rm}\CC[[u]]$ if $m\leq 0$.

ii) for all $m\geq 0$, $\alpha, \beta$ the coefficient of $\phi_{i,-m}^-.v_{\alpha}^{\lambda}$ on $v_{\beta}^{\lambda}$ is in $u^{-mR}\CC[[u]]$

iii) for all $m\geq 0$, $\alpha,\beta$ the coefficient of $\phi_{i,m}^+.v_{\alpha}^{\lambda}$ on $v_{\beta}^{\lambda}$ is in $\CC[[u]]$\end{defi}

\noindent Example : For $V\in\mathcal{O}(\U_q(\hat{\Glie}))$, we have $i(V)\in\mathcal{O}^0(\tilde{\U}_q(\hat{\Glie}))$.

\begin{thm}\label{prod} Let $V_1\in\mathcal{O}(\tilde{\U}_q(\hat{\Glie}))$ and $V_2\in\mathcal{O}^R(\tilde{\U}_q'(\hat{\Glie}))$. Then $\Delta_u$ defines a structure of $\tilde{\U}_q'(\hat{\Glie})$-module on $i(V_1)\otimes_{\mathcal{C}} V_2$ which is in $\mathcal{O}^{R+1}(\tilde{\U}_q'(\hat{\Glie}))$. Moreover the $l,u$-weights of $i(V_1)\otimes_{\mathcal{C}} V_2$ are of the form $(\lambda_1+\lambda_2,\gamma_1(z)\gamma_2(uz))$ where $(\lambda_1,\gamma_1)$ is a $l$-weight of $V_1$ and $(\lambda_2,\gamma_2)$ is a $l,u$-weight of $V_2$.\end{thm}

\noindent Remark : $\gamma(u)(z)=\gamma_1(z)\gamma_2(uz)$ means that for $i\in I, m\geq 0$ : 
$$\gamma_{i,\pm m}^{\pm}(u)=\underset{0\leq l\leq m}{\sum} (\gamma_1)_{i,\pm l}(u)(\gamma_2)_{i,\pm (m-l)}(u)u^{\pm (m-l)}$$

\demo As the definition of $\Delta_u$ involves infinite sums, we have to prove that the action formally defined by $\Delta_u$ makes sense on $V_1'\otimes_{\mathcal{C}} V_2$ where we denote $V_1'=i(V_1)$. Indeed the weight spaces of $V_1'$ and $V_2$ are finite dimensional and for $\lambda,\mu\in\Hlie^*$ we can use a $\CC$-base $(v_{\alpha}^{1,\lambda})_{\alpha}$ of $(V_1)_{\lambda}$ as a $\mathcal{C}$-base of $(V_1')_{\lambda}$ and the $\mathcal{C}$-basis $(v_{\alpha'}^{2,\mu})$ of $(V_2)_{\mu}$ given by the definition of $\mathcal{O}^R(\tilde{\U}_q'(\hat{\Glie}))$. So consider $\lambda,\mu\in\Hlie^*$, $i\in I$ and let us investigate the coefficients ($r\in\ZZ,m\geq 0$):

\noindent we have $x_{i,r}^+.((V_1')_{\lambda}\otimes (V_2)_{\mu})\subset (V_1')_{\lambda+\alpha_i}\otimes (V_2)_{\mu}\oplus (V_1')_{\lambda}\otimes (V_2)_{\mu +\alpha_i}$. 

on $(V_1')_{\lambda}\otimes (V_2)_{\mu +\alpha_i}$ : the coefficient of $x_{i,m}^+.(v_{\alpha}^{1,\lambda}\otimes v_{\alpha'}^{2,\mu})$ on $v_{\beta}^{1,\lambda}\otimes v_{\beta '}^{2,\mu+\alpha_i}$ is in $\underset{l\geq 0}{\sum}u^{r+l}\CC[[u]]\subset \CC[[u]]$ if $r\geq 0$, in $\underset{l\geq 0}{\sum}u^{r+l}u^{R(r+l)}\CC[[u]]\subset u^{(R+1)r}\CC[[u]]$ if $r\leq 0$.

on $(V_1')_{\lambda+\alpha_i}\otimes (V_2)_{\mu}$ : the coefficient of $x_{i,r}^+.(v_{\alpha}^{1,\lambda}\otimes v_{\alpha'}^{2,\lambda})$ on $v_{\beta}^{1,\lambda+\alpha_i}\otimes v_{\beta'}^{2,\mu}$ is in $\CC$.

\noindent we have $x_{i,r}^-.((V_1')_{\lambda}\otimes (V_2)_{\mu})\subset (V_1')_{\lambda-\alpha_i}\otimes (V_2)_{\mu}\oplus (V_1')_{\lambda}\otimes (V_2)_{\mu -\alpha_i}$. 

on $(V_1')_{\lambda}\otimes (V_2)_{\mu -\alpha_i}$ : the coefficient of $x_{i,r}^-.(v_{\alpha}^{1,\lambda}\otimes v_{\alpha'}^{2,\mu})$ on $v_{\beta}^{1,\lambda}\otimes v_{\beta'}^{2,\mu-\alpha_i}$ is in $u^r\CC[[u]]\subset \CC[[u]]$ if $r\geq 0$, in $u^ru^{Rr}\CC[[u]]$ if $r\leq 0$.

on $(V_1')_{\lambda-\alpha_i}\otimes (V_2)_{\mu}$ : the coefficient of $x_{i,r}^-.(v_{\alpha}^{1,\lambda}\otimes v_{\alpha'}^{2,\mu})$ on $v_{\beta}^{1,\lambda-\alpha_i}\otimes v_{\beta''}^{2,\mu}$ is in $\underset{l\geq 0}{\sum}u^l\CC[[u]]\subset \CC[[u]]$. 

\noindent we have $\phi_{i,\pm m}^{\pm}.((V_1')_{\lambda}\otimes (V_2)_{\mu})\subset ((V_1')_{\lambda}\otimes (V_2)_{\mu})$.

the coefficient of $\phi_{i,m}^+.(v_{\alpha}^{1,\lambda}\otimes v_{\alpha'}^{2,\mu})$ on $v_{\beta}^{1,\lambda}\otimes v_{\beta'}^{2,\mu}$ is in $\underset{0\leq l\leq m}{\sum}u^l \CC[[u]]\subset\CC[[u]]$.

the coefficient of $\phi_{i,-m}^-.(v_{\alpha}^{1,\lambda}\otimes v_{\alpha'}^{2,\mu})$ on $v_{\beta}^{1,\lambda}\otimes v_{\beta'}^{2,\mu}$ is in $\underset{0\leq l\leq m}{\sum}u^{-l}u^{-lR}\CC[[u]]\subset u^{-m(R+1)}\CC[[u]]$.

\noindent So we have a structure of $\tilde{\U}_q'(\hat{\Glie})$-module on $V_1'\otimes_{\mathcal{C}} V_2$. Let us prove that it is in $\mathcal{O}(\tilde{\U}_q'(\hat{\Glie}))$. We verify the properties of definition \ref{oprime}: i) ii) iii) are clear because the restriction of $\Delta_u$ to $\U_q(\hat{\Hlie})$ is the restriction of $\Delta_{\hat{\Hlie}}$. For iv) we note that for $(\lambda_1, \gamma_1), (\lambda_2,\gamma_2)\in P_{l,u}$, the $(V_1')_{\lambda_1, \gamma_1}\otimes (V_2)_{\lambda_2, \gamma_2}$ is in the pseudo weight space of $l,u$-weight $(\lambda_1+\lambda_2,\gamma_1(z)\gamma_2(zu))$ because $\Delta_u(\phi_i^{\pm}(z))=\phi_i^{\pm}(z)\otimes\phi_i^{\pm}(zu)$ (it also proves the last point of the proposition). 

\noindent Finally we see in the above computations that the coefficients verify the property of 
\\$\mathcal{O}^{R+1}(\tilde{\U}_q'(\hat{\Glie}))$, so $V_1'\otimes_{\mathcal{C}} V_2$ is in $\mathcal{O}^{R+1}(\tilde{\U}_q'(\hat{\Glie}))$. \qed

\begin{defi} For $R\geq 0$, we denote $\otimes_R : \mathcal{O}(\U_q(\hat{\Glie}))\times \mathcal{O}^R(\tilde{\U}_q'(\hat{\Glie}))\rightarrow\mathcal{O}^{R+1}(\tilde{\U}_q'(\hat{\Glie}))$ the bilinear map constructed in theorem \ref{prod}.\end{defi}

\noindent See section \ref{ex} for explicit examples. For $R\geq 2$ and $V_1, V_2,..., V_R \in\mathcal{O}(\U_q(\hat{\Glie}))$, one can define the iterated tensor product $V_1\otimes_{R-2} (V_2\otimes_{R-3}(...\otimes_0 V_R))...)$ which is in $\mathcal{O}^{R-1}(\tilde{\U}_q'(\hat{\Glie}))$.
 
\subsection{Simple modules of $\tilde{\U}_q'(\hat{\Glie})$} 

\subsubsection{$l,u$-highest weight modules} For $(\lambda,\Psi(u))\in P_{l,u}$, let $\tilde{M}(\lambda,\Psi(u))$ be the Verma $\tilde{\U}_q'(\hat{\Glie})$ module of highest weight $(\lambda,\Psi(u))$ (it is non trivial thanks to the triangular decomposition of $\tilde{\U}_q(\hat{\Glie})$ in lemma \ref{stepdeux}). So we have an analog of proposition \ref {simple} : for $(\lambda, \Psi(u))\in P_{l,u}$, there is a unique up to isomorphism simple $\tilde{\U}_q'(\hat{\Glie})$-module $\tilde{L}(\lambda, \Psi(u))$ of $l,u$-highest weight $(\lambda, \Psi(u))$ that is to say that there is $v\in \tilde{L}(\lambda, \Psi(u))$ such that ($i\in I, r\in\ZZ, m\geq 0, h\in \Hlie$):
$$x_{i,r}^+.v=0\text{ , }\tilde{L}(\lambda,\Psi(u))=\tilde{\U}_q'(\hat{\Glie}).v\text{ , }\phi_{i,\pm m}^{\pm}.v=\Psi_{i,\pm m}^{\pm}(u)v\text{ , }k_h.v=q^{\lambda(h)}.v$$
In a similar way one define the simple $\U_q'(\hat{\Glie})$-module $L(\lambda, \Psi(u))$ of $l,u$-highest weight $(\lambda, \Psi(u))$ (it is non trivial thanks to theorem \ref{dectrian}).

\begin{lem}\label{coree} For $(\lambda,\Psi(u))\in P_{l,u}$ we have an isomorphism of $\U_q(\hat{\Hlie})$-modules $\tilde{L}(\lambda,\Psi(u))\simeq L(\lambda,\Psi(u))$.\end{lem}

\demo Let $\tilde{M}'(\lambda,\Psi(u))\subset \tilde{M}(\lambda,\Psi(u))$ be the maximal proper $\tilde{\U}_q'(\hat{\Glie})$-submodule of $\tilde{M}'(\lambda,\Psi(u))$. It suffices to prove that $\tilde{\tau}_-.1$ is included in $\tilde{M}'(\lambda,\Psi(u))$ (see section \ref{finpreuve}; it implies that $\tilde{L}(\lambda,\Psi(u))$ is also a $\U_q'(\hat{\Glie})$-modules). It is a consequence of lemma \ref{steptrois}.\qed

\noindent In particular $\tilde{L}(\lambda,\Psi(u))\in\mathcal{O}(\tilde{\U}_q'(\hat{\Glie}))\Leftrightarrow L(\lambda,\Psi(u))\in\mathcal{O}(\U_q'(\hat{\Glie}))$ and in this case $\text{ch}_{q,u}(\tilde{L}(\lambda,\Psi(u)))=\text{ch}_{q,u}(L(\lambda,\Psi(u)))$.

\subsubsection{The category $\mathcal{O}_{\text{int}}(\tilde{\U}_q'(\hat{\Glie}))$} 

\begin{defi} $QP_{l,u}^+$ is the set of $(\lambda, \Psi(u))\in P_{l,u}$ satifying the following conditions :

i) for $i\in I$ there exist polynomials $Q_{i,u}(z)=(1-za_{i,1}u^{b_{i,1}})...(1-za_{i,N_i}u^{b_{i,N_i}}), R_{i,u}(z)=(1-zc_{i,1}u^{d_{i,1}})...(1-zc_{i,N_i'}u^{d_{i,N_i'}})$ ($a_{i,j}, c_{i,j}\in\CC^*$, $b_{i,j},d_{i,j}\geq 0$) such that in $\CC[u,u^{-1}][[z]]$ (resp. in 
\\$\CC[u,u^{-1}][[z^{-1}]]$):
$$\underset{r\geq 0}{\sum}\Psi_{i,\pm r}^{\pm}(u) z^{\pm r}= q_i^{\text{deg}(Q_{i,u})-\text{deg}(R_{i,u})}\frac{Q_{i,u}(zq_i^{-1})R_{i,u}(zq_i)}{Q_{i,u}(zq_i)R_{i,u}(zq_i^{-1})}$$

ii) there exist $\omega\in P^+$, $\alpha\in Q^+$ satisfying  $\lambda = \omega - \alpha$.

\noindent $P_{l,u}^+$ is the set of $(\lambda, \Psi(u))\in QP_{l,u}^+$ such that one can choose $R_{i,u}=1$ (in this case we denote $P_{i,u}=Q_{i,u}$). \end{defi}

\begin{lem} If $(\lambda ,\Psi(u))\in P_{l,u}^+$ then $\tilde{L}(\lambda,\Psi(u))\in\mathcal{O}(\tilde{\U}_q'(\hat{\Glie}))$. Moreover for $(\mu,\gamma(u))\in P_{l,u}$ we have $\text{dim}(\tilde{L}(\lambda,\Psi(u))_{\mu,\gamma(u)})\neq 0\Rightarrow (\mu,\gamma(u))\in QP_{l,u}^+$.\end{lem}

\noindent Remark : it follows from lemma \ref{coree} that we have the same results for $L(\lambda,\Psi(u))\in\mathcal{O}(\U_q'(\hat{\Glie}))$.

\demo Let $(\lambda ,\Psi(u))\in P_{l,u}^+$ and decompose $P_{i,u}(z)$ in the form :
$$P_{i,u}(z)=P_i^{(0)}(z)P_i^{(1)}(uz)...P_i^{(R)}(u^Rz)$$
where $R\geq 0$, $P_i^{(k)}(z)\in\CC[z]$, $P_i^{(k)}(0)=1$ for $0\leq k\leq R$ ($R$ can be taken large enough so that we have this form for all $i\in I$). For $0\leq k\leq R$, set $\Psi^{(k)}_i(z)=q_i^{\text{deg}(P_i^{(k)})}\frac{P_i^{(k)}(zq_i^{-1})}{P_i(zq_i)}$. For $1\leq k\leq R$ define $\lambda_k=\underset{i\in I}{\sum}\text{deg}(P_i^{(k)})\Lambda_i\in\Hlie^*$. Set $\lambda_0=\lambda-\underset{k=1...R}{\sum}\lambda_k$. Then for $0\leq k\leq R$ the $(\lambda_k,\Psi^{(k)})\in P_l^+$ and we can consider $L(\lambda_k,\Psi^{(k)})\in\mathcal{O}(\U_q(\hat{\Glie}))$. Let $V\in\mathcal{O}^R(\tilde{\U}_q'(\hat{\Glie}))$ be defined by :
$$V=i(L(\lambda_0,\Psi^{(0)}))\otimes_{R-1}(i(L(\lambda_1,\Psi^{(1)}))\otimes_{R-2}...(i(L(\lambda_{R-1},\Psi^{(R-1)}))\otimes_0 i(L(\lambda_R,\Psi^{(R)})))...)$$
Consider the $\tilde{\U}_q'(\hat{\Glie})$ submodule $L$ of $V$ generated by the tensor product of the highest weight vectors. It is a highest weight module of highest weight $(\lambda,\Psi(u))$. So $\tilde{L}(\lambda,\Psi(u))$ is a quotient of $L$ and so is in $\mathcal{O}(\tilde{\U}_q'(\hat{\Glie}))$.

\noindent For the second point it follows from proposition \ref{fr} that the $l,u$-weight of $i(L(\lambda_{k},\Psi^{(k)}))$ are in $QP_{l,u}^+$. So with the help of the last point of theorem \ref{prod} we see that the $l,u$-weights of $V$ are in $QP_{l,u}^+$ and we have the property for $\tilde{L}(\lambda,\Psi(u))$.\qed

\begin{defi} Let $\mathcal{O}_{\text{int}}(\tilde{\U}_q'(\hat{\Glie}))$ be the subcategory of modules $V\in\mathcal{O}(\tilde{\U}_q'(\hat{\Glie}))$ whose $l,u$-weights are in $QP_{l,u}^+$.\end{defi}

\begin{lem}\label{helpprop} For a module $V\in\mathcal{O}_{\text{int}}(\tilde{\U}_q'(\hat{\Glie}))$ there are $P_{(\lambda,\Psi(u))}\geq 0$ ($(\lambda,\Psi(u))\in P_{l,u}^+$) such that:
$$\text{ch}_{q,u}(V)=\underset{(\lambda,\Psi(u))\in P_{l,u}^+}{\sum}P_{(\lambda,\Psi(u))}\text{ch}_{q,u}(\tilde{L}(\lambda,\Psi(u)))=\underset{(\lambda,\Psi(u))\in P_{l,u}^+}{\sum}P_{(\lambda,\Psi(u))}\text{ch}_{q,u}(L(\lambda,\Psi(u)))$$\end{lem}

\demo Analogous to the proof of proposition \ref{helpthm} (the second identity follows from lemma \ref{coree}).\qed

\subsection{$\CC[u^{\pm}]$-forms and specialization}

\subsubsection{$\CC[u^{\pm}]$-forms}

Let $\U_q^u(\hat{\Glie})=\U_q(\hat{\Glie})\otimes_{\CC}\CC[u^{\pm}]\subset\U_q'(\hat{\Glie})$.

\begin{defi} A $\CC[u^{\pm}]$-form of a $\U_q'(\hat{\Glie})$-module $V$ is a sub-$\U_q^u(\hat{\Glie})$-module $L$ of $V$ such that the map $\mathcal{C}\otimes_{\CC[u^{\pm}]}L\rightarrow V$ is an isomorphism of $\U_q'(\hat{\Glie})$-module.\end{defi}

\noindent Note that it means that $L$ generates $V$ as $\mathcal{C}$-vector space and that some vectors which are $\CC[u^{\pm}]$-linearly independent in $L$ are $\mathcal{C}$-linearly independent in $V$.

\noindent Let us look at some examples:

\begin{prop}\label{intform} For $(\lambda,\Psi(u))\in P_{l,u}$ and $v$ a highest weight vector of the Verma module $M(\lambda, \Psi(u))$ (resp. the simple module $L(\lambda,\Psi(u))$), the $\U_q^u(\hat{\Glie})$-module $\U_q^u(\hat{\Glie}).v$ is a $\CC[u^{\pm}]$-form of $M(\lambda,\Psi(u))$ (resp. of $L(\lambda,\Psi(u))$) which is isomorphic to the Verma (resp. the simple) $\U_q^u(\hat{\Glie})$-module of l,u-highest weight $(\lambda,\Psi(u))$.\end{prop}

\demo As $(\lambda,\Psi(u))$ is fixed, we omit it. $M$ is the quotient of $\mathcal{C}\otimes_{\CC}\U_q(\hat{\Glie})$ by the relations generated by $x_{i,r}^{\pm}=\phi_{i,\pm m}^{\pm}-\Psi_{i,\pm m}^{\pm}(u)=k_h -q^{\lambda(h)}=0$. So the relations between monomials are in $\CC[u^{\pm}]$ and $\U_q^u(\hat{\Glie}).1\subset M$ is a $\CC[u^{\pm}]$-form of $M$. Moreover those relations are the same as in the construction of the Verma $\U_q^u(\hat{\Glie})$-module $M^u$ as a quotient of $\CC[u^{\pm}]\otimes_{\CC}\U_q(\hat{\Glie})$; and so $\U_q^u(\hat{\Glie}).1\simeq M^u$.

\noindent Let us look at $L$. Denote by $L^u$ the simple $\U_q^u(\hat{\Glie})$-module of highest weight $(\lambda,\Psi(u))$. We have $L=M/M'$ (resp. $L^u=M^u/M'^u$) where $M'$ (resp. $M'^u$) in the maximal proper submodule of $M$ (resp. $M^u$). 

\noindent The $\mathcal{C}$-subspace $M''$ of $M$ generated by $M'^u$ is
isomorphic to $\mathcal{C}\otimes_{\CC[u^{\pm}]}M'^u$ (because $M^u$ is a
$\CC[u^{\pm}]$-form of $M$). As $M''$ has no vector of weight $\lambda$,
it is a proper submodule of $M$ and $M''\subset M'$. Suppose that $M'\neq
M''$ and consider $M'/M''\subset M/M''$. $M^u/M'^u$ is a
$\CC[u^{\pm}]$-form of $M/M''$. Let $v$ be a non zero highest weight
vector of $M'/M''$ and let us write:
$v=\underset{\alpha}{\sum}f_{\alpha}(u)v_{\alpha}$ where $v_{\alpha}\in
M^u/M'^u$ and $f_{\alpha}(u)\in\mathcal{C}$ (as there is a finite number 
of $f_{\alpha}(u)$, we can suppose that they are $\CC[u^{\pm}]$-linearly 
independent). For all $i\in I,r\in\ZZ$, we have $x_{i,r}^+.v=0$ and so
for all $\alpha$, $x_{i,r}^+.v_{\alpha}=0$. Fix $w_{\alpha}\in M^u$ whose image in $M^u/M'^u$ is $v_{\alpha}$. As for all $i\in I, r\in\ZZ$, $x_{i,r}^+.w_{\alpha}\in M'^u$, $\U_q^u(\hat{\Glie}).w_{\alpha}$ is a proper submodule of $M^u$ and $w_{\alpha}\in M'^u$. So $v=0$, contradiction. So $M'=M''$. In particular $M'\simeq M'^u\otimes_{\CC[u^{\pm}]}\mathcal{C}$, $M'\cap M^u=M'^u$.

\noindent For $v$ a highest weight vector of $L$, the $\U_q^u(\hat{\Glie}).v\simeq \U_q^u(\hat{\Glie}).1=M^u/(M^u\cap M')=M^u/M'^u=L^u$ is a $\CC[u^{\pm}]$-form of $L$.\qed

\subsubsection{Specializations}\label{spe} Consider $p:\mathcal{E}_{l,u}\rightarrow \mathcal{E}_l$ the surjection such that $p((\lambda,\Psi(u)))=(\lambda,\Psi(1))$.

\begin{lem}\label{sper} Let $V$ be in $\mathcal{O}(\U_q'(\hat{\Glie}))$. If $L$ is a $\CC[u^{\pm}]$-form of $V$ then the specialization $L'=L/(1-u)L$ of $L$ is in $\mathcal{O}(\U_q(\hat{\Glie}))$ and $\text{ch}_q(L')=p(\text{ch}_{q,u}(V))$.\end{lem}

\demo Indeed for $(\mu,\gamma(u))\in QP_{l,u}$ consider $L_{\mu,\gamma(u)}=L\cap V_{\mu,\gamma(u)}$. As $p:L\otimes_{\CC[u]}\mathcal{C}\rightarrow V$ is an isomorphism, we have $V_{\mu,\gamma(u)}\simeq p^{-1}(V_{\mu,\gamma(u)})=L_{\mu,\gamma(u)}\otimes_{\CC[u]}\mathcal{C}$. In particular $L_{\mu,\gamma(u)}$ is a free $\CC[u^{\pm}]$ of rank $\text{dim}_{\mathcal{C}}(V_{\mu,\gamma(u)})$. So $\text{dim}_{\CC}(L'_{\mu})=\text{dim}_{\mathcal{C}}(V_{\mu})$, and $L'\in\mathcal{O}(\U_q(\hat{\Glie}))$. We can conclude because:
$$L_{\lambda,\gamma}'=\underset{(\lambda,\gamma(u))\in p^{-1}((\lambda,\gamma))}{\bigoplus}(L_{\lambda,\gamma(u)}/(u-1)L_{\lambda,\gamma(u)})$$\qed

\subsubsection{Proof of theorem \ref{posc}} For $(\lambda,\Psi(u))\in P_{l,u}^+$, it follows from proposition \ref{intform} and lemma \ref{sper} that $p(\text{ch}_{q,u}(L(\lambda,\Psi(u))))$ is of the form $\text{ch}_q(L)$ where $L\in\mathcal{O}_{\text{int}}(\U_q(\hat{\Glie}))$, that is to say $p(\text{ch}_{q,u}(L(\lambda,\Psi(u))))\in\text{ch}_q(\text{Rep}^+(\U_q(\hat{\Glie})))$. So (lemma \ref{helpprop}) for $V\in\mathcal{O}_{\text{int}}(\tilde{U}_q'(\hat{\Glie}))$ we have $p(\text{ch}_{q,u}(V))\in\text{ch}_q(\text{Rep}^+(\U_q(\hat{\Glie})))$. 

\noindent Consider $V_1,V_2\in\mathcal{O}_{\text{int}}(\U_q(\hat{\Glie}))$. We have seen that $p(\text{ch}_{q,u}(i(V_1)\otimes_0 i(V_2))\in\text{ch}_q(\text{Rep}^+(\U_q(\hat{\Glie})))$. But : 
$$p(\text{ch}_{q,u}(i(V_1)\otimes_0 i(V_2))=\text{ch}_q(V_1)\text{ch}_q(V_2)$$ 
because the specialization of $\Delta_u$ on $\U_q(\hat{\Hlie})$ at $u=1$ is $\Delta_{\hat{\Hlie}}$. This ends the proof of theorem \ref{posc}.\qed

\subsection{Example}\label{ex} We study in detail an example in the case $\Glie=sl_2$ where everything is computable thanks to Jimbo's evaluation morphism (see \cite{Cha, Cha2}). In this case we have $\U_q(\hat{sl_2})=\tilde{\U}_q(\hat{sl_2})$.

\noindent For $a\in\CC^*$ consider $V=L(1-za)\in\mathcal{O}_{\text{int}}(\U_q(\hat{sl_2}))$. $V$ is two dimensional $V=\CC v_0\oplus \CC v_1$ and for $r\in \ZZ$, $m\geq 1$ the action of $\U_q(\hat{sl_2})$ is given in the following table:

$$
\begin{array}{l|l|l}
   & v_0 & v_1\\ \hline
 x_r^+ & 0 & a^r v_0 \\ \hline
 x_r^- & a^r v_1 & 0 \\\hline 
 \phi_{\pm m }^{\pm} & \pm (q-q^{-1}) a^{\pm m} v_0 & \mp (q-q^{-1}) a^{\pm m} v_1 \\ \hline

k^{\pm} & q^{\pm} v_0 & q^{\mp} v_1 \\ \hline
\phi^{\pm}(z) & q\frac{1-q^{-2}az}{1-az} v_0 & q^{-1}\frac{1-q^2az}{1-az}v_1
\\ 
\end{array}  
$$
\noindent Remark : in the table $\phi^{\pm}(z)\in\U_q(\hat{\Glie})[[z^{\pm}]]$ acts on $V[[z^{\pm}]]$.

\noindent For $a,b\in\CC^*$, let $V=L(1-za), W=L(1-zb)\in \mathcal{O}_{\text{int}}(\U_q(\hat{\Glie}))$. Consider basis $V=\CC v_0\oplus\CC v_1$, $W=\CC w_0\oplus \CC w_1$ as in the previous table. The tensor product $\otimes_0$ defines an action of $\U_q'(\hat{sl_2})$ on $X=i(V)\otimes_{\mathcal{C}} i(W)$ (see theorem \ref{prod}). $X$ is a 4 dimensional $\mathcal{C}$-vector space of base $\{v_0\otimes w_0, v_1\otimes w_0, v_0\otimes w_1, v_1\otimes w_1\}$. The action of $\U_q'(\hat{sl_2})$ is given by ($r\in\ZZ$):

$$
\begin{array}{l|l|l}
   & v_0\otimes w_0 & v_1\otimes w_0 \\ \hline
 x_r^+ & 0 & a^r (v_0\otimes w_0) \\ \hline
x_{r}^- & u^rb^r(v_0\otimes w_1)+a^rq\frac{1-q^{-2}a^{-1}ub}{1-a^{-1}ub}(v_1\otimes w_0)& u^rb^r(v_1\otimes w_1)\\ \hline
\phi^{\pm}(z) & q^2\frac{(1-q^{-2}az)(1-q^{-2}buz)}{(1-az)(1-buz)}(v_0\otimes w_0)& \frac{(1-q^2 az)(1-q^{-2}buz)}{(1-az)(1-buz)}(v_1\otimes w_0)
\end{array}  
$$

$$
\begin{array}{l|l|l}
   & v_0\otimes w_1 & v_1\otimes w_1\\ \hline
 x_r^+ & b^ru^r \frac{1-q^2a^{-1}bu}{1-uba^{-1}}(v_0\otimes w_0) & a^r(v_0\otimes w_1)+b^ru^rq\frac{1-q^{-2}a^{-1}bu}{1-uba^{-1}}(v_1\otimes w_0)\\ \hline
x_{r}^- & a^r q^{-1}\frac{1-q^2ua^{-1}b}{1-a^{-1}ub}(v_1\otimes w_1)& 0 \\ \hline
\phi^{\pm}(z) & \frac{(1-q^{-2} az)(1-q^2buz)}{(1-az)(1-buz)}(v_0\otimes w_1)& q^{-2}\frac{(1-q^2 az)(1-q^2 buz)}{(1-az)(1-buz)}(v_1\otimes w_1)
\\ 
\end{array}  
$$
\noindent Remark : in the table $\phi^{\pm}(z)\in\U_q(\hat{\Glie})[[z^{\pm}]]$ acts on $X[[z^{\pm}]]$.

\noindent Consider the $l$-weights $\gamma_a,\gamma_a',\gamma_b,\gamma_b'\in P_l$ (the $\lambda\in\Hlie^*$ can be omitted because $sl_2$ is finite) :
$$\gamma_a^{\pm}(z)=q\frac{1-q^{-2}az}{1-az}\text{ , }\gamma_a'^{\pm}(z)=q^{-1}\frac{1-q^2az}{1-az}\text{ , }\gamma_b^{\pm}(z)=q\frac{1-q^{-2}bz}{1-bz}\text{ , }\gamma_b'^{\pm}(z)=q^{-1}\frac{1-q^2bz}{1-bz}$$
Consider also $\gamma_a(z)\gamma_b(uz),\gamma_a'(z)\gamma_b(uz),\gamma_a(z)\gamma_b'(uz),\gamma_a'(z)\gamma_b'(uz)\in P_{l,u}$. We see that :
$$\text{ch}_{q,u}(X)=e(\gamma_a(z)\gamma_b(uz))+e(\gamma_a'(z)\gamma_b(uz))+e(\gamma_a(z)\gamma_b'(uz))+e(\gamma_a(z)\gamma_b'(uz))$$
Those $l,u$-weights are distinct, the $l,u$-weights spaces are 1 dimensional:
$$X=(X)_{\gamma_a(z)\gamma_b(uz)}\oplus(X)_{\gamma_a'(z)\gamma_b(uz)}\oplus(X)_{\gamma_a(z)\gamma_b'(uz)}\oplus(X)_{\gamma_a'(z)\gamma_b'(uz)}$$
We see that $X$ is of highest weight $\gamma_a(z)\gamma_b(uz)\in P_{l,u}$. Let us prove that it is simple : indeed $X$ has no proper submodule : if for all $r\in \ZZ$, $x_r^+.(\alpha (v_1\otimes w_0)+\beta (v_0\otimes w_1))=0$, then for all $r\in\ZZ$, $\alpha a^r+\beta b^ru^r \frac{1-q^2a^{-1}bu}{1-uba^{-1}}=0$. In particular $\alpha +\beta \frac{1-q^2a^{-1}bu}{1-uba^{-1}}=0$ and $a^r-b^ru^r=0$ for all $r\in\ZZ$, impossible. So $X\simeq L(\gamma_a(z)\gamma_b(uz))$ as a $\U_q'(\hat{sl_2})$-module. It follows from proposition \ref{intform} that $\tilde{X}=\U_q^u(\hat{\Glie}).(v_0\otimes w_0)\subset X$ is a $\CC[u^{\pm}]$-form of $X$.

\noindent Let us look explicitly at this $\CC[u^{\pm}]$-form : consider $e_1,e_2,e_3,e_4\in \tilde{X}$ defined by: 
$$e_1=v_0\otimes w_0\text{ , }e_2=x_0^-.e_1\text{ , }e_3=-a^{-1}x_1^-.e_1+e_2\text{ , }e_4=qx_0^-.e_2$$
We have the following formulas:
$$e_1=v_0\otimes w_0\text{ , }e_2=(v_0\otimes w_1)+q\frac{1-q^{-2}a^{-1}bu}{1-a^{-1}ub}(v_1\otimes w_0)\text{ , }e_3=(1-uba^{-1})(v_0\otimes w_1)\text{ , }e_4=(v_1\otimes w_1)$$
Moreover the action of $\U_q^u(\hat{\Glie})$ is given by ($r\in\ZZ$):
$$
\begin{array}{l|l|l}
   & e_1 & e_2 \\ \hline
 x_r^+ & 0 & (qa^r\frac{1-(a^{-1}ub)^{r+1}}{1-a^{-1}ub}-q^{-1}bua^{r-1}\frac{1-(a^{-1}ub)^{r-1}}{1-a^{-1}ub})e_1 \\ \hline
x_{r}^- & a^r(e_2-\frac{1-(a^{-1}ub)^r}{1-a^{-1}ub}e_3)& (q^{-1}a^r\frac{1-(a^{-1}ub)^{r+1}}{1-a^{-1}ub}-q bua^{r-1}\frac{1-(a^{-1}ub)^{r-1}}{1-a^{-1}ub})e_4\\ \hline
\phi^{\pm}(z) & q^2\frac{(1-q^{-2}az)(1-q^{-2}buz)}{(1-az)(1-buz)}e_1 & \frac{(1-q^2 az)(1-q^{-2}buz)}{(1-az)(1-buz)}e_2+\frac{az(q^2-q^{-2})}{(1-az)(1-buz)}e_3  \\
\end{array}  
$$

$$
\begin{array}{l|l|l}
   & e_3 & e_4 \\ \hline
 x_r^+ & b^ru^rq^{-1}(1-q^2a^{-1}bu)e_1 &b^ru^r e_2+a^r\frac{1-(uba^{-1})^r}{1-uba^{-1}}e_3 \\ \hline
x_{r}^- & a^rq^{-1}(1-q^2ua^{-1}b)e_4 & 0 \\ \hline
\phi^{\pm}(z) & \frac{(1-q^{-2}az)(1-q^2 buz)}{(1-az)(1-buz)}e_3 & q^{-2}\frac{(1-q^2 az)(1-q^2 buz)}{(1-az)(1-buz)}e_4
\\ 
\end{array}  
$$
In particular we see that $\CC[u^{\pm}]e_1\oplus\CC[u^{\pm}]e_2\oplus\CC[u^{\pm}]e_3\oplus\CC[u^{\pm}]e_4$ is stable by the action of $\U_q^u(\hat{\Glie})$, so is equal to $\tilde{X}$. So we have verified that $X\simeq \tilde{X}\otimes_{\CC[u^{\pm}]}\mathcal{C}$.

\noindent Let us describe the specialization of $\tilde{X}$ at $u=1$ : let $\tilde{X}'=\CC e_1\oplus\CC e_2\oplus \CC e_3\oplus \CC e_4$. The action of $\U_q(\hat{\Glie})$ on $\tilde{X}'$ is given by (for $z\in\CC, r\in\ZZ$, we denote $[z]_r'=\frac{1-z^r}{1-z}\in\ZZ[z^{\pm}]$ ($z\neq 1$) and $[1]_r'=r$):

$$
\begin{array}{l|l|l}
   & e_1 & e_2 \\ \hline
 x_r^+ & 0 & (qa^r[a^{-1}b]_{r+1}'-q^{-1}ba^{r-1}[a^{-1}b]_{r-1}')e_1 \\ \hline
x_{r}^- & a^r(e_2-[a^{-1}b]_r' e_3)& (q^{-1}a^r[a^{-1}b]_{r+1}'-q ba^{r-1}[a^{-1}b]_{r-1}')e_4\\ \hline
\phi^{\pm}(z) & q^2\frac{(1-q^{-2}az)(1-q^{-2}bz)}{(1-az)(1-bz)}e_1 & \frac{(1-q^2 az)(1-q^{-2}bz)}{(1-az)(1-bz)}e_2+\frac{az(q^2-q^{-2})}{(1-az)(1-bz)}e_3  \\
\end{array}  
$$

$$
\begin{array}{l|l|l}
   & e_3 & e_4 \\ \hline
 x_r^+ & b^rq^{-1}(1-q^2a^{-1}b)e_1 &b^r e_2+a^r[a^{-1}b]_r'e_3 \\ \hline
x_{r}^- & a^mq^{-1}(1-q^2a^{-1}b)e_4 & 0 \\ \hline
\phi^{\pm}(z) & \frac{(1-q^{-2}az)(1-q^2 bz)}{(1-az)(1-bz)}e_3 & q^{-2}\frac{(1-q^2 az)(1-q^2 bz)}{(1-az)(1-bz)}e_4
\\ 
\end{array}  
$$
We see that $\tilde{X}'=\U_q(\hat{\Glie}).e_1$. Moreover if $ab^{-1}\notin \{q^2,q^{-2}\}$ : $\tilde{X}'$ has no proper submodule because the formula $x_m^+(\alpha e_2+\beta e_3)=0$ means that for all $r\in\ZZ$:
$$\alpha(qa^r[a^{-1}b]_{r+1}-q^{-1}ba^{r-1}[a^{-1}b]_r)+\beta b^rq^{-1}(1-q^2a^{-1}b)=0$$
which is possible only if $ab^{-1}\in\{q^2,q^{-2}\}$ or $\alpha=\beta=0$. So:

if  $ab^{-1}\notin\{q^2,q^{-2}\}$, $\tilde{X}'\simeq L(\gamma_a\gamma_b)$ is simple and :
$$\text{ch}_q(V)\text{ch}_q(W)=\text{ch}_q(\tilde{X}')=\text{ch}_q(L(\gamma_a\gamma_b))$$

if $ab^{-1}=q^2$ (resp. $ab^{-1}=q^{-2}$), $\CC e_3\subset \tilde{X}'$ (resp. $\CC((q^2-1)e_2+e_3)\subset \tilde{X}'$) is a submodule of $\tilde{X}'$ isomorphic to $L(1)$ and :
$$\text{ch}_q(V)\text{ch}_q(W)=\text{ch}_q(\tilde{X}')=\text{ch}_q(L(\gamma_a\gamma_b))+\text{ch}_q(L(1))$$

\end{document}